\documentclass[11pt]{article}
\usepackage{amsfonts}
\usepackage{bbm}
\usepackage{amscd}
\usepackage{xypic}
\usepackage{amsmath}
\usepackage{amssymb}
\usepackage{amsthm}
\usepackage{bbm}
\usepackage{CJK}
\usepackage{fancyhdr}
\usepackage{hyperref}
\usepackage{indentfirst}
\usepackage{latexsym}
\usepackage{mathrsfs}
 \allowdisplaybreaks 
\usepackage{color}

\usepackage[top=1in,bottom=1in,left=1.25in,right=1.25in]{geometry}
\def\<{\langle}
\def\>{\rangle}
\def\a{\alpha}
\def\b{\beta}
\def\ci{\circ}
\def\c{\cdot}

\def\D{\Delta}

\def\i{\iota}

\def\r{\rho}
\def\lr{\longrightarrow}

\def\o{\otimes}

\def\vp{\varphi}

\def\<{\langle}
\def\>{\rangle}

\date{}
\begin{document}
\renewcommand{\baselinestretch}{1.2}
\renewcommand{\arraystretch}{1.0}
\title{\bf Constructing New Braided $T$-Categories via Weak Monoidal Hom-Hopf Algebras}
\date{}
\author{{\bf Wei Wang,
        Shuanhong Wang\footnote {Corresponding author:Shuanhong Wang, E-mail: shuanhwang2002@yahoo.com}
        and Xiaohui Zhang }\\
{\small Department of Mathematics, Southeast University}\\
{\small Jiangsu Nanjing 210096, P. R. CHINA}\\}
 \maketitle
\begin{center}
\begin{minipage}{14.cm}
\begin{center}{\bf ABSTRACT}\end{center}

In this paper, we define and study weak monoidal Hom-Hopf algebras, which generalize both weak Hopf algebras and monoidal Hom-Hopf algebras.
If $H$ is a weak monoidal Hom-Hopf algebra with bijective antipode and let $Aut_{wmHH}(H)$ be the set of all automorphisms of $H$. Then we introduce a category ${_{H}\mathcal{WMHYD}^{H}}(\alpha,\beta)$ with $\alpha,\beta\in Aut_{wmHH}(H)$ and construct a braided $T$-category $\mathcal{WMHYD}(H)$ that having all the categories ${_{H}\mathcal{WMHYD}^{H}}(\alpha,\beta)$ as components.

 \vskip 0.5cm

\end{minipage}
\end{center}

\section*{ Introduction}
\def\theequation{0. \arabic{equation}}
\setcounter{equation} {0} \hskip\parindent

Hom-algebras first appeared in \cite{AS081}, where the associativity was replaced by the Hom-associativity and similar to the Hom-coassociativity (see in \cite{AS10, AS082}). Based on these properties, definitions of Hom-bialgebras, Hom-Hopf algebras and further developments existed later in \cite{G10}, \cite{AF14}-\cite{AS082}, \cite{D09}, \cite{D10}, \cite{SI11} and \cite{FGS09}. In \cite{SI11}, the authors illustrated Hom-structures from the point of view of monoidal categories and introduced monoidal Hom-algebras, monoidal Hom-coalgebras, etc., in a symmetric monoidal category, which were different from the Hom-algebras and Hom-coalgebras.\\

Weak Hom-Hopf algebra was introduced in \cite{ZW}. The axioms were the same as the ones for a Hom-Hopf
algebra, except that the coproduct of the unit, the product of counit and the antipode condition were replaced by the weaker properties.\\

Turaev in \cite{V94,V08} generalized quantum invariants if 3-manifolds to the case of a 3-manifold $M$ endowed with a homotopy class of maps
$M\rightarrow K(G,1)$, where $G$ is a group. And braided $T$-categories which are braided monoidal categories in Freyd-Yetter categories of crossed $G$-sets(see in \cite{PD89}) play a key role of constructing these homotopy invariants.\\

Then one main question naturally arises. How to construct classes of new braided $T$-categories? \\

The main propose of this paper is to generalize the notion of monoidal Hom-Hopf algebras to weak cases. We first define the notion of weak monoidal Hom-bigebras and weak momoidal Hom-Hopf algebras. And we give examples to show its not a trivial generalization of monoidal Hom-Hopf algebras.
We also investigate the relationship between weak monoidal Hom-Hopf algebras and weak Hom-Hopf algebras.
After that, we construct new examples of braided $T$-categories, which generalize the construction given by Panaite and Staic (see in \cite{PM07}).\\

The article is organized as follows.\\

In Section 1, we recall definitions and basic results related to monoidal Hom-Hopf algebras and braided $T$-categories.
In Section 2, we apply the construction  given by Stef and Isar \cite{SI11} to the category of vector spaces,
then we give definitions of weak monoidal Hom-bialgebras and some axioms which will be used in Section 3 and Section 4. We also define the weak monoidal Hom-Hopf algebras and give some properties of its antipode.\\

In Section 3, we introduce a class of new categories $ _{H}\mathcal{WMHYD}^{H}(\alpha, \beta)$ (see Definition 3.1) of weak $(\alpha,\beta)$-Yetter-Drinfeld modules associated with $\alpha,\beta\in Aut_{wmHH}(H)$. Furthermore, we prove that the category $ _{H}\mathcal{WMHYD}^{H}(\alpha, \beta)$ is actually a weak monoidal entwined Hom-module category $_{H}\mathcal{M}^{H}(\psi(\alpha, \beta))$.
Then in Section 4,  we prove $_{H}\mathcal{WMHYD}^{H}$ is a monoidal cateogory and then construct a class of new braided $T$-categories $_{H}\mathcal{WMHYD}^{H}$ in the sense of Turaev\cite{V08}.\\

\section*{1. Preliminaries}
\def\theequation{1. \arabic{equation}}
\setcounter{equation} {0} \hskip\parindent

Throughout, let $k$ be a fixed field. Everything is over $k$ unless
 otherwise specified.  We
 refer the readers to the books of Sweedler \cite{M69}
  for the relevant concepts on the general theory of Hopf
 algebras.  Let $(C, \Delta )$ be a coalgebra. We use the Sweedler-Heyneman's notation for
 $\Delta $ as follows:
 \begin{eqnarray*}
 \Delta (c)=\sum c_1\otimes c_2,
 \end{eqnarray*}
 for all $c\in C$.

  \vskip 0.5cm
{\bf 1.1. Monoidal Hom-Hopf algebras.}
\vskip 0.5cm

Let $\mathcal{M}_{k}=(\mathcal{M}_{k},\o,k,a,l,r )$
 denote the usual monoidal category of $k$-vector spaces and linear maps between them.
 Recall from \cite{SI11}
 that there is the {\it monoidal Hom-category} $\widetilde{\mathcal{H}}(\mathcal{M}_{k})=
 (\mathcal{H}(\mathcal{M}_{k}),\,\o,\,(k,\,id),
 \,\widetilde{a},\,\widetilde{l},\,\widetilde{r })$, a new monoidal category,
  associated with $\mathcal {M}_{k}$ as follows:

 $\bullet$  The objects of
 $ \mathcal{H}(\mathcal{M}_{k})$ are couples
 $(M,\mu)$, where $M \in \mathcal {M}_{k}$ and $\mu \in Aut_k(M)$, the set of
   all $k$-linear automomorphisms of $M$;

  $\bullet$  The morphism $f:(M,\mu)\rightarrow (N,\nu)$ in $ \mathcal{H}(\mathcal{M}_{k})$
  is the $k$-linear map $f: M\rightarrow N$ in $\mathcal{M}_{k}$
  satisfying   $ \nu \circ f = f\ci \mu$, for any two objects
  $(M,\mu),(N,\nu)\in \mathcal{H}(\mathcal{M}_{k})$;

 $\bullet$   The tensor product is given by
 $$
 (M,\mu)\o (N,\nu)=(M\o N,\mu\o\nu )
 $$
for any $(M,\mu),(N,\nu)\in \mathcal{H}(\mathcal{M}_{k})$.

$\bullet$ The tensor unit is given by $(k, id)$;

$\bullet$   The associativity constraint $\widetilde{a}$
 is given by the formula
 $$
 \widetilde{a}_{M,N,L}=a_{M,N,L}\circ((\mu\o id)\o
  \varsigma^{-1})=(\mu\o(id\o\varsigma^{-1}))\circ a_{M,N,L},
  $$
 for any objects $(M,\mu),(N,\nu),(L,\varsigma)\in \mathcal{H}(\mathcal{M}_{k})$;

 $\bullet$  The left and right unit constraint
  $\widetilde{l}$ and $\widetilde{r }$ are given by
 $$
 \widetilde{l}_M=\mu\circ l_M=l_M\circ(id\o\mu),\, \quad
  \widetilde{r}_M =\mu\circ r_M=r_M\circ(\mu\o id)
  $$
for all $(M,\mu) \in \mathcal{H}(\mathcal{M}_{k})$.
\\

We  now recall from \cite{SI11} the following notions used later.
\\

 A {\it unital monoidal Hom-associative algebra} (a monoidal Hom-algebra in
 Proposition 2.1 of \cite{SI11}) is a vector space $A$
 together with an element $1_A\in A$ and linear maps
$$m:A\o A\rightarrow A;\,\,a\o b\mapsto ab, \,\,\,\alpha\in Aut_k(A)$$
such that
\begin{equation}
\alpha(a)(bc)=(ab)\alpha(c),
  \end{equation}
$$\alpha(ab)=\alpha(a)\alpha(b),$$
  \begin{equation}
  a1_A=1_Aa=\alpha(a),
  \end{equation}
\begin{equation}
\alpha(1_A)=1_A ,
  \end{equation}
for all $a,b,c\in A.$
\\

{\bf Remark.}   (1) In the language of Hopf algebras, $m$ is called the
 Hom-multiplication, $\alpha$ is the twisting
 automorphism and $1_A$ is the unit. Note that Eq.(1.1) can be rewirtten as
 $a(b\alpha^{-1}(c)) = (\alpha^{-1}(a)b)c$.
  The monoidal Hom-algebra $A$ with $\alpha$ will be denoted by $(A,\alpha)$.

 (2)  Let $(A,\alpha)$ and  $(A',\alpha')$ be two monoidal Hom-algebras.
   A monoidal Hom-algebra map
 $f:(A,\alpha)\rightarrow (A',\alpha')$ is a linear map such that
 $f\circ \alpha=\alpha'\circ f,f(ab)=f(a)f(b)$ and
 $f(1_A)=1_{A'}.$

(3) The definition of monoidal Hom-algebras is different
 from the unital Hom-associative algebras in \cite{AS10} and
 \cite{AS082} in the following sense.
  The same twisted associativity condition (1.1) holds in both cases.
 However, the unitality condition in their notion is the usual untwisted one:
  $a1_A=1_Aa =a,$ for any $a\in A,$ and the twisting map $\alpha$
  does not need to be monoidal (that is, (1.2) and (1.3) are not required).
\\

  {\it A counital monoidal Hom-coassociative coalgebra} is
  an object $(C,\gamma)$ in the category
  $\tilde{\mathcal{H}}(\mathcal{M}_{k})$
 together with linear maps
 $\D:C\rightarrow C\o C,\,\D(c)=c_1\o c_2$ and
 $\varepsilon:C\rightarrow k$ such that
 \begin{equation}
\gamma^{-1}(c_1)\o\D(c_2)=\D(c_1)\o\gamma^{-1}(c_2),
  \end{equation}
 \begin{equation}
\D(\gamma(c))=\gamma(c_1)\o\gamma(c_2),
  \end{equation}
 $$c_1\varepsilon(c_2)=\gamma^{-1}(c)=\varepsilon(c_1)c_2,$$
 \begin{equation}
 \varepsilon(\gamma(c))=\varepsilon(c)
  \end{equation}
 for all $c\in C.$
\\

 {\bf Remark.}   (1) Note that (1.4)is equivalent to
 $c_1\o c_{21}\o \gamma(c_{22})=\gamma(c_{11})\o c_{12}\o c_2.$
 Analogue to monoidal Hom-algebras, monoidal Hom-coalgebras
 will be short for counital monoidal Hom-coassociative coalgebras
 without any confusion.

 (2)  Let $(C,\gamma)$ and $(C',\gamma')$ be two monoidal Hom-coalgebras.
 A monoidal Hom-coalgebra map $f:(C,\gamma)\rightarrow(C',\gamma')$
 is a linear map such that $f\circ \gamma=\gamma'\circ f, \D\circ f=(f\o f)\circ\D$
  and $\varepsilon'\circ f=\varepsilon.$
\\

 {\it A monoidal Hom-bialgebra} $H=(H,\alpha,m,1_H,\D,\varepsilon)$
 is a bialgebra in the monoidal category
 $ \tilde{\mathcal{H}}(\mathcal {M}_{k}).$
 This means that $(H,\alpha,m,1_H)$ is a monoidal Hom-algebra and
 $(H,\alpha,\D,\varepsilon)$ is
 a monoidal Hom-coalgebra such that $\D$ and $\varepsilon$
  are morphisms of algebras,
 that is, for all $h,g\in H,$
 $$\D(hg)=\D(h)\D(g),\, \,\, \D(1_H)=1_H\o1_H,\,\,\,\,\,\,
  \varepsilon(hg)=\varepsilon(h)\varepsilon(g), \,\,\,\,\,\varepsilon(1_H)=1.$$\\

  A monoidal Hom-bialgebra $(H,\alpha)$ is called {\it a monoidal Hom-Hopf algebra}
 if there exists a morphism (called antipode)
 $S: H\rightarrow H$ in $ \tilde{\mathcal{H}}(\mathcal {M}_{k})$
 (i.e., $S\ci \alpha=\alpha\ci S$),
 which is the convolution inverse of the identity morphism $id_H$
 (i.e., $ S*id=1_H\ci \varepsilon=id*S$). Explicitly,  for all $h\in H$,
 $$
 S(h_1)h_2=\varepsilon(h)1_H=h_1S(h_2).
 $$\\

 {\bf Remark.}   (1) Note that a monoidal Hom-Hopf algebra is
 by definition a Hopf algebra in $ \tilde{\mathcal{H}}(\mathcal {M}_{k})$.

 (2)  Furthermore, the antipode of monoidal Hom-Hopf algebras has
   almost all the properties of antipode of Hopf algebras such as
 $$S(hg)=S(g)S(h),\,\,\,\, S(1_H)=1_H,\,\,\,\,
  \D(S(h))=S(h_2)\o S(h_1),\,\,\,\,\,\,\varepsilon\ci S=\varepsilon.$$
 That is, $S$ is a monoidal Hom-anti-(co)algebra homomorphism.
  Since $\alpha$ is bijective and commutes with $S$,
  we can also have that the inverse $\alpha^{-1}$ commutes with $S$,
  that is, $S\ci \alpha^{-1}= \alpha^{-1}\ci S.$
 \\

 In the following, we recall the notions of actions on monoidal
 Hom-algebras and coactions on monoidal Hom-coalgebras.
\\

  Let $(A,\alpha)$ be a monoidal Hom-algebra.
{\it A left $(A,\alpha)$-Hom-module} consists of
 an object $(M,\mu)$ in $\tilde{\mathcal{H}}(\mathcal {M}_{k})$
  together with a morphism
  $\psi:A\o M\rightarrow M,\psi(a\o m)=a\cdot m$ such that
 $$\alpha(a)\c(b\c m)=(ab)\c\mu(m),\,\,
 \,\,\mu(a\c m)=\alpha(a)\c\mu(m),\,\,
 \,\,1_A\c m=\mu(m),$$
 for all $a,b\in A$ and $m \in M.$
\\

 Monoidal Hom-algebra $(A,\alpha)$ can be
 considered as a Hom-module on itself by the Hom-multiplication.
  Let $(M,\mu)$ and $(N,\nu)$ be two left $(A,\alpha)$-Hom-modules.
  A morphism $f:M\rightarrow N$ is called left
   $(A,\alpha)$-linear if
   $f(a\c m)=a\c f(m),f\ci \mu= \nu\ci f.$
  We denoted the category of left $(A,\alpha)$-Hom modules by
  $\tilde{\mathcal{H}}(_{A}\mathcal {M}_{k})$.
\\

 Similarly, let $(C,\gamma)$ be a monoidal Hom-coalgebra.
 {\it A right $(C,\gamma)$-Hom-comodule} is an object
  $(M,\mu)$ in $\tilde{\mathcal{H}}(\mathcal {M}_{k})$
  together with a $k$-linear map
  $\rho_M:M\rightarrow M\o C,\rho_M(m)=m_{(0)}\o m_{(1)}$ such that
 \begin{equation}
 \mu^{-1}(m_{(0)})\o \D_C(m_{(1)})
 =(m_{(0)(0)}\o m_{(0)(1)})\o \gamma^{-1}(m_{(1)}),
  \end{equation}
    \begin{equation}
 \rho_M(\mu(m))=\mu(m_{(0)})\o\gamma(m_{(1)}),
  \ \ \
 m_{(0)}\varepsilon(m_{(1)})=\mu^{-1}(m),
 \end{equation}
for all $m\in M.$
\\

 $(C,\gamma)$ is a Hom-comodule on itself via the Hom-comultiplication.
 Let $(M,\mu)$ and $(N,\nu)$ be two right $(C,\gamma)$-Hom-comodules.
  A morphism $g:M\rightarrow N$ is called right $(C,\gamma)$-colinear
  if $g\ci \mu=\nu\ci g$ and
  $g(m_{(0)})\o m_{(1)}=g(m)_{(0)}\o g(m)_{(1)}.$
  The category of right
  $(C,\gamma)$-Hom-comodules is denoted by
   $\tilde{\cal{H}}(\cal {M}^C)$ .
\\

 Let $(H, \alpha)$ be a monoidal Hom-bialgebra. We now recall from \cite{YL14}
 that a monoidal Hom-algebra $(B,\beta)$ is called {\it a left
 $H$-Hom-module algebra}, if $(B,\beta)$ is a left
 $H$-Hom-module with action $\cdot$ obeying the following axioms:
\begin{equation}
h\c(ab)=(h_1\c a)(h_2\c b),\,\,\,\,\,\,\,\,h\c1_B=\varepsilon(h)1_B,
  \end{equation}
for all $a,b\in B,h\in H.$
\\

 Recall from \cite{LB14} that a monoidal Hom-algebra $(B,\beta)$ is called
 {\it a left $H$-Hom-comodule algebra}, if $(B,\beta)$ is a left
 $H$-Hom-comodule with coaction $\rho$ obeying the following axioms:
 $$\rho(ab) = a_{(-1)}b_{(-1)}\otimes a_{(0)}b_{(0)},
  \ \ \ \rho_l(1_B)= 1_B\otimes 1_H,$$
for all $a,b\in B,h\in H.$
\\

 Let $(H,m,\Delta,\alpha)$ be a monoidal Hom-bialgebra.
 Recall from (\cite{LB14, YW14}) that a {\it left-right Yetter-Drinfeld Hom-module}
 over $(H,\alpha)$ is the object $(M,\cdot,\rho,\mu)$ which
 is both in  $\tilde{\cal{H}}(_{H}\cal {M})$ and $\tilde{\cal{H}}(\cal {M}^{H})$
 obeying the compatibility condition:
 \begin{equation}
 h_{1}\c m_{(0)}\o h_{2}m_{(1)}=(\alpha(h_{2})\c m)_{(0)} \o
 \alpha^{-1}(\alpha(h_{2})\c m)_{(1)})h_{1}.
  \end{equation}

 {\bf Remark.} (1) The category of all left-right Yetter-Drinfeld Hom-modules
  is denoted by $\tilde{\cal{H}}(_{H}\cal {YD} ^{H})$ with understanding morphism.

 (2) If $(H,\alpha)$ is a monoidal Hom-Hopf algebra with a bijective
 antipode $S$, then the above equality is equivalent to
 $$\r (h\c m)=\alpha(h_{21})\c m_{(0)}\o (h_{22}\a^{-1}(m_{(1)}))S^{-1}(h_{1}), $$
 for all $h\in H$ and $m\in M$.\\

 \vskip 0.5cm
 {\bf 1.2. Braided $T$-categories.}
\vskip 0.5cm

 A {\sl monoidal category} $\mathcal{C}=(\mathcal{C},\mathbb{I},\otimes,a,l,r)$
 is a category $\mathcal{C}$ endowed with a functor
 $\otimes: \mathcal{C}\times\mathcal{C}\rightarrow\mathcal{C}$
 (the {\sl tensor product}), an object $\mathbb{I}\in \mathcal{C}$
 (the {\sl tensor unit}), and natural isomorphisms $a$
 (the {\sl associativity constraint}), where
 $a_{U,V,W}:(U\otimes V)\otimes W\rightarrow U\otimes (V\otimes W)$
 for all $U,V,W\in \mathcal{C}$, and $l$ (the {\sl left unit constraint})
 where $l_U: \mathbb{I}\otimes U\rightarrow U,\,r$
 (the {\sl right unit constraint}) where
 $r_{U}:U\otimes\mathcal{C}\rightarrow U$ for all $U\in \mathcal{C}$,
 such that for all $U,V,W,X\in \mathcal{C},$
 the {\sl associativity pentagon}
 $a_{U,V,W\otimes X}\circ a_{U\otimes V,W,X}
 =(U\otimes a_{V,W,X})\circ a_{U,V\otimes W,X}\circ
 (a_{U,V,W}\otimes X)$ and
 $(U\otimes l_V)\circ(r_U\otimes V)=a_{U,I,V}$ are satisfied.
 A monoidal categoey $\mathcal{C}$ is {\sl strict} when all
  the constraints are identities.
\\

  Let $G$ be a group and let $Aut(\mathcal{C})$ be the group of
  invertible strict tensor functors from $\mathcal{C}$ to itself.
  A category $\mathcal{C}$ over $G$
  is called a {\sl crossed category } if it satisfies the following:

 $\bullet$  $\mathcal{C}$ is a monoidal category;

 $\bullet$  $\mathcal{C}$ is disjoint union of a family
 of subcategories $\{\mathcal{C}_{\alpha }\}_{\alpha \in G}$, and for any $U\in \mathcal{C}_{\alpha }$,
 $V\in \mathcal{C}_{\beta }$, $U\otimes V\in \mathcal{C}_{\a \beta }$.
 The subcategory $\mathcal{C}_{\alpha}$ is called the $\alpha$th component of $\mathcal{C}$;

 $\bullet$  Consider a group homomorphism $\varphi : G\rightarrow Aut(\mathcal{C})$, $\b \mapsto \varphi _{\b }$, and
  assume that$\varphi_{\beta }(\varphi_{\alpha })=\varphi_{\beta\alpha\beta^{-1}}$,for all $\alpha,\beta\in G$.
  The functors $\varphi_{\beta }$  are called conjugation isomorphisms.

 Furthermore, $\mathcal{C}$ is called strict when it is
 strict as a monoidal category.
\\

 {\sl Left index notation}: Given $\a \in G$
 and an object $V\in \mathcal{C}_{\a }$, the functor $\vp _{\a }$
 will be denoted by ${}^V( \cdot )$, as in Turaev \cite{V08} or
 Zunino \cite{M04}, or even ${}^{\a }( \cdot )$.
 We use the notation ${}^{\overline{V}}( \cdot )$
 for ${}^{\a ^{-1}}( \cdot )$. Then we have
 ${}^V id_U=id_{{} V^U}$ and
 ${}^V(g\circ f)={}^Vg\circ {}^Vf$.
 Since the conjugation $\vp : G\lr Aut(\mathcal{C})$ is a
 group homomorphism, for all $V, W\in \mathcal{C}$, we have ${}^{V\o W}( \cdot )
 ={}^V({}^W( \cdot ))$ and ${}^\mathbb{I}( \cdot )={}^V({}^{\overline{V}}( \cdot ))
 ={}^{\overline{V}}({}^V( \cdot ))=id_\mathcal{C}$. Since, for all
 $V\in \mathcal{C}$, the functor ${}^V( \cdot )$ is strict, we have
 ${}^V(f\o g)={}^Vf\o {}^Vg$, for any morphisms $f$ and $g$ in $\mathcal{C}$,
  and ${}^V\mathbb{I}=\mathbb{I}$.
\\

 A {\sl braiding} of a crossed category $\mathcal{C}$ is
 a family of isomorphisms $({c=c_{U,V}})_{U,V}\in \mathcal{C}$,
 where $c_{U,V}: U\otimes V\rightarrow {}^UV\otimes U$
 satisfying the following conditions:\\
 (i)~For any arrow $f\in \mathcal{C}_{\a }(U, U')$ and
 $g\in \mathcal{C}(V, V')$,
 $$
 (({}^{\a }g)\o f)\circ c _{U, V}=c _{U' V'}\circ (f\o g).
 $$
 (ii)~For all $ U, V, W\in \mathcal{C},$ we have
 $$
 c _{U\o V, W}=a_{{}^{U\o V}W, U, V}\circ (c _{U, {}^VW}\o
 id_V)\circ a^{-1}_{U, {}^VW, V}\circ (\i _U\o c _{V, W})
  \circ a_{U, V, W},
  $$
    $$c _{U, V\o W}=a^{-1}_{{}^UV, {}^UW, U}
 \circ (\i _{({}^UV)}\o c _{U, W})\circ a_{{}^UV, U, W}\ci
 (c _{U, V}\o \i_W)\circ a^{-1}_{U, V, W},$$
 where $a$ is the natural isomorphisms in the tensor category
 $\mathcal{C}$.\\
 (iii)~For all $ U, V\in \mathcal{C}$ and $\b\in G$,
 $$ \vp _{\b }(c
 _{U, V})=c _{\vp _{\b }(U), \vp _{\b }(V)}.
 $$

 A crossed category endowed with a braiding is called
 a {\sl braided $T$-category}.\\

\section*{2. Weak monoidal Hom-bialgebras and examples}
\def\theequation{2. \arabic{equation}}
\setcounter{equation} {0} \hskip\parindent

Let $k$ be a commutative ring. The results of Hom-construction by Stef and Isar\cite{SI11} can be applied to the category of $k$-modules (vector spaces if $k$ is a field) $\mathcal{C}=\mathcal{M}_{k}$.

In this section we will introduce the notion of a weak monoidal Hom-bialgebra.
\\

{\bf Definition 2.1.}~ $H=(H,\xi,m, 1_H,\Delta,\varepsilon)$ is called a {\sl weak monoidal Hom-bialgebra} if $(H,\xi)$ is both a monoidal Hom-algebra and a monoidal Hom-coalgebra, satisfying the following identities for any $a,b,c \in H$:
\[\left.\begin{array}{l}
(1)~~\Delta(ab)  = \Delta(a)\Delta(b);\\
(2)~~\varepsilon((ab)c) = \varepsilon(ab_1)\varepsilon(b_2c),~~\varepsilon(a(bc)) = \varepsilon(ab_2)\varepsilon(b_1c);\\
(3)~~(\Delta \o id_H)\Delta(1_H) = 1_1 \o 1_21'_{1} \o 1'_{2},~~(id_H \o \Delta)\Delta(1_H) = 1_1 \o1'_{1}1_2 \o 1'_{2}.
\end{array}\right.
\]

{\bf Definition 2.2.}~ If $(H,\xi)$ and $(H',\xi')$ are two weak monoidal Hom-bialgebras, a linear map $f:H\rightarrow H'$ is called a {\sl morphism of weak monoidal Hom-bialgebras} if $f$ is both a morphism of monoidal Hom-algebras and a morphism of monoidal Hom-coalgebras.
\\

Let $H$ be a weak monoidal Hom-bialgebra. Define linear maps $\varepsilon_s$ and $\varepsilon_t$ by the formulas
$$\aligned
\varepsilon_s(h) = \xi^{2}(1_{1})\varepsilon(\xi^{-2}(h)1_{2}),\quad \varepsilon_t(h) = \varepsilon(1_{1}\xi^{-2}(h))\xi^{2}(1_{2}),
\endaligned$$
for any $h \in H$, where $\varepsilon_t$, $\varepsilon_s$ are called the \textsl{target} and \textsl{source counital maps}. We adopt the notations $H_t = \varepsilon_t(H)$ and $H_s = \varepsilon_s(H)$ for their images.

Similarly, we define the linear maps $\widehat{\varepsilon_{s}}$ and $\widehat{\varepsilon_{t}} $ by the formulas
$$\aligned
\widehat{\varepsilon_s}(h) = \xi^{2}(1_{1})\varepsilon(1_{2}\xi^{-2}(h)),\quad \widehat{\varepsilon_t}(h) = \varepsilon(\xi^{-2}(h) 1_{1})\xi^{2}(1_{2}),
\endaligned$$
for any $h \in H$. Their images are denoted by $\widehat{H_t} = \widehat{\varepsilon_t}(H)$ and $\widehat{H_s} = \widehat{\varepsilon_s}(H)$.
And then we obtain the following identities
$$
\Delta(1_H) \in H_s \o H_t,~~\Delta(1_H) \in \widehat{H_s} \o \widehat{H_t},
$$
and
$$1_1 \o 1_2 = \xi(1_1) \o \xi(1_2),~~\xi(h_1) \o \xi(h_2) = 1_1h_1 \o 1_2h_2= h_11_1 \o h_21_2.
$$\\

{\bf Proposition 2.3.}~ Let $H$ be a weak monoidal Hom-bialgebra. Then for any $h \in H$, we have
\begin{eqnarray}
&&\mbox{(i)}~~~~\varepsilon_s(h_1) \o h_2 = \xi^{3}(1_1) \o \xi^{-2}(h)1_2;\\
&&\mbox{(ii)}~~~\widehat{\varepsilon_s}(h_1) \o h_2 = \xi(1_1) \o 1_2\xi^{-2}(h).
\end{eqnarray}

{\bf Proof. }
(i). From the definition of $\varepsilon_s$, we immediately get that
$$\aligned
\varepsilon_s(h_1) \o h_2 &= \xi^{2}(1_1)\varepsilon(\xi^{-2}(h_1)1_{2}) \o h_2\\
&=\xi^{2}(1_{1})\varepsilon(\xi^{-2}(h_{1})(1'_{1}\xi^{-1}(1_{2})))\otimes \xi^{-1}(h_{2})1'_{2}\\
&=\xi^{2}(1_{1})\otimes \xi^{-2}(h)\varepsilon(1'_{1}1_{2})1'_{2}\\
&=\xi^{2}(1_{1})\otimes \xi^{-2}(h)\xi^{-1}(1_{2})\\
&=\xi^{3}(1_{1})\otimes \xi^{-2}(h)1_{2}.
\endaligned$$

The proof of (ii) is similar to (i).
 $\hfill \blacksquare$
 \\

{\bf Theorem 2.4.}~ Let $H$ be a weak Hom-bialgebra. Then for any $a,b,c \in H$, we have the following identities
\begin{eqnarray}
&&\mbox{(i)}~~~\Delta(1_1) \o 1_2 = 1_1 \o \Delta(1_2);\\
&&\mbox{(ii)}~~\varepsilon((ab)c) = \varepsilon(a(bc)).
\end{eqnarray}
{\bf Proof. }
From the Proposition 2.3 above, we know that $\varepsilon_{s}(1_{1})\otimes 1_{2}=\xi^{3}(1_{1})\otimes \xi(1_{2})$, another side,
$$\aligned
\varepsilon_s(1_1) \o 1_2 &= \xi^{2}(1'_1)\varepsilon(1_{1}1'_2) \o \xi^{2}(1_{2})\\
&=\xi^{2}(1_{1})\otimes \xi(1_{2})
\endaligned$$
then we can get
\begin{eqnarray}
1_1 \otimes 1_2=\xi (1_1) \otimes 1_2.
\end{eqnarray}

Similarly, we can get $\widehat{\varepsilon_s}(1_1) \o 1_2 = \xi(1_1) \o \xi(1_2)$ from (2.2) and $\widehat{\varepsilon_s}(1_1) \o 1_2 = \xi(1_1) \o \xi^{2}(1_2) $ by a direct computation, which means
\begin{eqnarray}
1_{1} \otimes 1_{2}=1_1 \otimes \xi(1_{2}).
\end{eqnarray}

(i). We compute as follows:
$$\aligned
\Delta(1_1) \o 1_2 \stackrel {(2.6)}{=} \Delta(1_1) \o \xi^{-1}(1_2) = \xi^{-1}(1_1) \o \Delta(1_2) \stackrel {(2.5)}{=} 1_1 \o \Delta(1_2).
\endaligned$$

(ii). Since
$$\aligned
\varepsilon(\xi(a)b) = \varepsilon(a1_1)\varepsilon(1_2b) = \varepsilon(a1_{1}))\varepsilon(\xi(1_{2})\xi(b)) \stackrel {(2.6)}{=}\varepsilon(a1_1)\varepsilon(1_2\xi(b)) = \varepsilon(\xi(a)\xi(b))=\varepsilon(ab),
\endaligned$$
we get
\begin{eqnarray}
\varepsilon(\xi(a)b) = \varepsilon(ab),
\end{eqnarray}
and similarly, we can get
\begin{eqnarray}
\varepsilon(a\xi(b)) = \varepsilon(ab).
\end{eqnarray}
Thus we have
$$\aligned
\varepsilon(a(bc)) \stackrel {(2.7)}{=}\varepsilon(\xi(a)(bc)) = \varepsilon((ab)\xi(c)) \stackrel {(2.8)}{=}\varepsilon((ab)c).
\endaligned$$
 $\hfill \blacksquare$
 \\

So from now on, we can denote $(\Delta\o id_{H})\Delta(1_H)$ by $1_{1} \o 1_{2} \o 1_{3}$, and $\varepsilon((ab)c)$ by $ \varepsilon(abc)$.\\

From (2.5) and (2.6), we can rewrite the definitions of $\varepsilon_{s}$, $\varepsilon_{t}$, $\widehat{\varepsilon_s}$ and $\widehat{\varepsilon_t}$ by the following formulas:
\begin{eqnarray}
&&\varepsilon_{s}(h)=1_{1}\varepsilon(h1_{2})\quad \varepsilon_{t}(h)=\varepsilon(1_{1}h)1_{2},\\
&&\widehat{\varepsilon_{s}}(h)=1_{1}\varepsilon(1_{2}h)\quad \widehat{\varepsilon_{t}}(h)=\varepsilon(h 1_{1})1_{2}
\end{eqnarray}\\

{\bf Corollary 2.5.}~ Let $H$ be a weak monoidal Hom-bialgebra. Then for any $a,b,c \in H$, we have
\begin{eqnarray}
&&\mbox{(i)}~~~~\varepsilon_s(ab) = \varepsilon_s(\xi(a)b) = \varepsilon_s(a\xi(b)),\quad \varepsilon_t(ab) = \varepsilon_t(\xi(a)b) = \varepsilon_t(a\xi(b));\\
&&~~~~~~~\widehat{\varepsilon_s}(ab) = \widehat{\varepsilon_s}(\xi(a)b) = \widehat{\varepsilon_s}(a\xi(b)),\quad\widehat{\varepsilon_t}(ab) = \widehat{\varepsilon_t}(\xi(a)b) = \widehat{\varepsilon_t}(a\xi(b));\\
&&\mbox{(ii)}~~~\varepsilon_s\circ \varepsilon_s = \varepsilon_s\circ\xi = \xi\circ \varepsilon_s = \varepsilon_s,\quad \varepsilon_t\circ \varepsilon_t = \varepsilon_t\circ\xi = \xi\circ \varepsilon_t = \varepsilon_t;\\
&&\mbox{(iii)}~~\widehat{\varepsilon_s}\circ \widehat{\varepsilon_s} = \widehat{\varepsilon_s}\circ\xi = \xi\circ \widehat{\varepsilon_s} = \widehat{\varepsilon_s},\quad \widehat{\varepsilon_t}\circ \widehat{\varepsilon_t} = \widehat{\varepsilon_t}\circ\xi = \xi\circ \widehat{\varepsilon_t} = \widehat{\varepsilon_t}.
\end{eqnarray}

{\bf Proof.} Using equations (2.5) to (2.10), left to readers.\\

{\bf Proposition 2.6.}
Let $H$ be a weak monoidal Hom-bialgebra.  Then for any $x,y,h \in H$,  we have
\begin{eqnarray}
&&\mbox{(i)}~~~~\varepsilon(xy) = \varepsilon(\varepsilon_s(x)y) = \varepsilon(x\varepsilon_t(y)) ;\\
&&~~~~~~~\varepsilon(xy) = \varepsilon(\widehat{\varepsilon_t}(x)y) = \varepsilon(x\widehat{\varepsilon_s}(y));\\
&&\mbox{(ii)}~~~x_{1} \o \varepsilon_t (x_{2})=1_{1}\xi^{-2}(x)\o1_{2},~~\varepsilon_s (x_1) \o x_2 = 1_1 \o \xi^{-2}(x)1_2;\\
&&~~~~~~~\widehat{\varepsilon_s}(x_{1}) \o x_{2}=1_{1}\o1_{2}\xi^{-2}(x),~~x_{1} \o \widehat{\varepsilon_t} (x_{2})=\xi^{-2}(x) 1_{1}\o1_{2};\\
&&\mbox{(iii)}~~x_{1}\varepsilon(h x_{2}) =\varepsilon_s (h)\xi^{-2}(x),~~\varepsilon(x_1h)x_2 = \xi^{-2}(x)\varepsilon_t(h) ;\\
&&~~~~~~~\varepsilon(h x_{1})x_{2} = \widehat{\varepsilon_t} (h)\xi^{-2}(x),~~x_1\varepsilon(x_2h) = \xi^{-2}(x)\widehat{\varepsilon_s}(h);\\
&&\mbox{(iv)}~~\varepsilon_s(x)\varepsilon_t(y) =\varepsilon_t(y) \varepsilon_s(x),~~\widehat{\varepsilon_s}(x)\widehat{\varepsilon_t}(y) = \widehat{\varepsilon_t}(y) \widehat{\varepsilon_s}(x).
\end{eqnarray}

{\bf Proof. } We just check some of them and the others are left to readers.

(i).$$\aligned
\varepsilon(\varepsilon_s(x)y) \stackrel {(2.9)}{=} \varepsilon(1_1y)\varepsilon(x 1_2) = \varepsilon(x\xi(y)) \stackrel {(2.8)}{=} \varepsilon(xy).
\endaligned$$

(ii). Similar to Proposition 2.3.

(iii).$$\aligned
x_1\varepsilon(hx_2) \stackrel {(2.15)}{=}x_1\varepsilon(h\varepsilon_t(x_2)) = 1_1\xi^{-2}(x)\varepsilon(h1_2) = \varepsilon_{s}(h)\xi^{-2}(x).
\endaligned$$

(iv).$$\aligned
\varepsilon_s(x)\varepsilon_t(y) = 1_1\varepsilon_t(y)\varepsilon(x1_2) \stackrel {(2.19)}{=} \varepsilon(1_1y)1_2\varepsilon(x1_3) = \varepsilon_t(y) \varepsilon_s(x).
\endaligned$$
 $\hfill \blacksquare$
 \\

{\bf Proposition 2.7.}
Let $H$ be a weak monoidal Hom-bialgebra.  Then for any $h \in H$,  we have
\begin{eqnarray}
&& \Delta (\varepsilon_t(h)) = 1_{1}\varepsilon_t (h) \o 1_{2} = \varepsilon_t(h) 1_{1} \o 1_{2};\\
&& \Delta (\varepsilon_s (h)) = 1_{1} \o 1_{2}\varepsilon_s (h) = 1_{1} \o\varepsilon_s (h) 1_{2};\\
&& \Delta (\widehat{\varepsilon_s} (h)) =1_{1} \o\widehat{\varepsilon_s} (h)1_{2} =1_{1} \o 1_{2}\widehat{\varepsilon_s} (h);\\
&& \Delta (\widehat{\varepsilon_t} (h))= 1_{1}\widehat{\varepsilon_t}(h) \o 1_{2} = \widehat{\varepsilon_t} (h)1_{1} \o 1_{2}.
\end{eqnarray}

{\bf Proof. } We only check Eq.(2.22). Indeed,
$$\aligned
\Delta(\varepsilon_t(h)) = \varepsilon(1_1h)1_2 \o 1_3 \stackrel {(2.19)}{=} 1_1\varepsilon_t(h) \o 1_2 \stackrel {(2.21)}{=}\varepsilon_t(h) 1_{1} \o 1_{2}.
\endaligned$$
The other three identities can be proved by similar calculations.
 $\hfill \blacksquare$\\

{\bf Proposition 2.8.}
Let $H$ be a weak monoidal Hom-bialgebra.  Then for any $x,y \in H$,  we have
\begin{eqnarray}
&&\mbox{(i)}~~~~\varepsilon_t(x\varepsilon_t(y)) = \varepsilon_t(xy),~~\varepsilon_t(\varepsilon_t(x)y) = \varepsilon_t(x)\varepsilon_t(y) ;\\
&&~~~~~~~\varepsilon_s(x\varepsilon_s(y)) = \varepsilon_s(x)\varepsilon_s(y),~~\varepsilon_s(\varepsilon_s(x)y) = \varepsilon_s(xy) ;\\
&&\mbox{(ii)}~~~\widehat{\varepsilon_t}(x\widehat{\varepsilon_t}(y)) = \widehat{\varepsilon_t}(x)\widehat{\varepsilon_t}(y), ~~\widehat{\varepsilon_t}(\widehat{\varepsilon_t}(x)y) =\widehat{\varepsilon_t}(xy) ;\\
&&~~~~~~~\widehat{\varepsilon_s}(x\widehat{\varepsilon_s}(y)) = \widehat{\varepsilon_s}(xy),~~\widehat{\varepsilon_s}(\widehat{\varepsilon_s}(x)y) = \widehat{\varepsilon_s}(x)\widehat{\varepsilon_s}(y);\\
&&\mbox{(iii)}~~\varepsilon_t(x\widehat{\varepsilon_s}(y)) = \varepsilon_t(xy),~~\varepsilon_t(\widehat{\varepsilon_t}(x)y)=\widehat{\varepsilon_t}(x)\varepsilon_t(y);\\
&&~~~~~~~\varepsilon_s(x\widehat{\varepsilon_s}(y)) = \varepsilon_s(x)\widehat{\varepsilon_s}(y),~~\varepsilon_s(\widehat{\varepsilon_t}(x)y) = \varepsilon_s(xy) ;\\
&&\mbox{(iv)}~~~\widehat{\varepsilon_t}(x\varepsilon_t(y)) = \widehat{\varepsilon_t}(x)\varepsilon_t(y),~~\widehat{\varepsilon_t}(\varepsilon_s(x)y) =
\widehat{\varepsilon_t}(xy) ;\\
&&~~~~~~~\widehat{\varepsilon_s}(x\varepsilon_t(y)) = \widehat{\varepsilon_s}(xy),~~\widehat{\varepsilon_s}(\varepsilon_s(x)y) =\varepsilon_s(x)  \widehat{\varepsilon_s}(y);\\
&&\mbox{(v)}~~~~\varepsilon_t (x_{1}) \o x_{2} = \varepsilon_t (1_{1}) \o 1_{2}\xi^{-2}(x),~~
  x_{1}\o \varepsilon_s (x_{2}) = \xi^{-2}(x) 1_{1}\o \varepsilon_s (1_{2});\\
&&~~~~~~~x_{1}\o \widehat{\varepsilon_s} (x_{2}) = 1_{1}\xi^{-2}(x) \o \widehat{\varepsilon_s} (1_{2}),~~\widehat{\varepsilon_t} (x_{1}) \o x_{2} = \widehat{\varepsilon_t} (1_{1}) \o \xi^{-2}(x) 1_{2}.
\end{eqnarray}

{\bf Proof. } Straightforward.
\\

{\bf Remark.}~ For any weak monoidal Hom-bialgebra $(H,\xi)$, if $\xi = id_H$, then $H$ is exactly a weak bialgebra. If $\Delta$ and $\varepsilon$ are all monoidal Hom-algebra morphisms, then $H$ is exactly a monoidal Hom-bialgebra.\\

Hereby we give examples of weak monoidal Hom-bialgebras.\\

{\bf Example 2.9.} If $(H,\xi,m,\eta,\Delta,\varepsilon)$ is a finite dimensional weak monoidal Hom-bialgebra and $H^\ast$ is the linear dual of $H$. Then
$(H^\ast,\xi^\ast,m^\ast,\eta^\ast,\Delta^\ast,\varepsilon^\ast)$ is also a finite dimensional weak monoidal Hom-bialgebra, where
$$
\xi^\ast(f) = f\circ\xi^{-1},
$$
for any $f \in H^\ast$.
\\

{\bf Example 2.10.} (10-dimensional weak monoidal Hom-bialgebra) Let $H$ be a vector space over $k$ with basis $\{x_{i}\}$, where $i=1,2,\cdots 10$ and let $0\neq\lambda\in k$, we give the structure of $H$ as follows:

$\bullet $ the multiplication

$$\begin{array}{|c|c|c|c|c|c|c|c|c|c|c|}
\hline H      &  x_1  & x_2  &x_3  &x_4  &x_5 & x_6  &x_7  &x_8  &x_9  &x_{10}  \\
\hline x_1    & x_1 & x_2 & x_3 & \lambda x_4 & \lambda x_5 & x_6  &x_7  &x_8  &\lambda x_9  &\lambda x_{10} \\
\hline x_2    & x_2 & x_2 & x_3 & \lambda x_4 & \lambda x_5 & x_7 & x_7 & x_8  &\lambda x_9  &\lambda x_{10}\\
\hline x_3    & x_3 & x_3 & x_2 & -\lambda x_5 & -\lambda x_4 & x_8 & x_8 & x_7 & -\lambda x_{10} & -\lambda x_9 \\
\hline x_4    & \lambda x_4 & \lambda x_4 & \lambda x_5 & 0 & 0 & \lambda x_9 & \lambda x_9 & \lambda x_{10} & 0 & 0\\
\hline x_5    & \lambda x_5 & \lambda x_5 & \lambda x_4 & 0 & 0 & \lambda x_{10} & \lambda x_{10} & \lambda x_9 & 0 & 0\\
\hline x_6    & x_6 & x_7 & x_8 & \lambda x_9 & \lambda x_{10} & x_6 & x_7 & x_8 & \lambda x_9 & \lambda x_{10}\\
\hline x_7    & x_7 & x_7 & x_8 & \lambda x_9 & \lambda x_{10} & x_7 & x_7 & x_8 & \lambda x_9 & \lambda x_{10}\\
\hline x_8    & x_8 & x_8 & x_7 & -\lambda x_{10} & -\lambda x_9 & x_8 & x_8 & x_7 & -\lambda x_{10} & -\lambda x_9\\
\hline x_9    & \lambda x_9 & \lambda x_9 & \lambda x_{10} & 0 & 0 & \lambda x_9 & \lambda x_9 & \lambda x_{10} & 0 & 0\\
\hline x_{10} & \lambda x_{10} & \lambda x_{10} & \lambda x_9 & 0 & 0 & \lambda x_{10} & \lambda x_{10} & \lambda x_9 & 0 & 0\\
\hline
 \end{array}$$

$\bullet $ the comultiplication

$\Delta(x_1)=x_1\otimes x_1-x_1\otimes x_6 -x_6\otimes x_1 +2x_6\otimes x_6 -x_1\otimes x_2 +x_1\otimes x_7 +x_6\otimes x_2$

$-2x_6\otimes x_7 -x_2\otimes x_1 +x_2\otimes x_6 +x_7\otimes x_1 -2x_7\otimes x_6+2x_2\otimes x_2 -2x_2\otimes x_7 -2x_7\otimes x_2 +4x_7\otimes x_7,$

$\Delta(x_2)= x_2\otimes x_2 -x_2\otimes x_7 -x_7\otimes x_2 +2x_7\otimes x_7,$

$\Delta(x_3)=x_3\otimes x_3 -x_3\otimes x_8 -x_8\otimes x_3+2x_8\otimes x_8,$

$\Delta(x_4)=\frac{1}{\lambda} (x_3\otimes x_4- x_3\otimes x_9 - x_8\otimes x_4 +2 x_8\otimes x_9 +  x_4\otimes x_2 - x_4\otimes x_7 - x_9\otimes x_2 +2x_9\otimes x_7),$

$\Delta(x_5)=\frac{1}{\lambda} (x_2\otimes x_5 - x_2\otimes x_{10} - x_7\otimes x_5 +2 x_7\otimes x_{10} + x_{5}\otimes x_{3} - x_5\otimes x_8 - x_{10}\otimes x_8 +2 x_{10}\otimes x_8),$

$\Delta(x_6)=x_6\otimes x_6-x_6\otimes x_7-x_7\otimes x_6+2x_7\otimes x_7$

$\Delta(x_7)=x_7\otimes x_7,$

$\Delta(x_8)=x_8\otimes x_8,$

$\Delta(x_9)=\frac{1}{\lambda} (x_8\otimes x_9+x_9\otimes x_7),$

$\Delta(x_{10})=\frac{1}{\lambda} (x_7\otimes x_{10}+x_{10}\otimes x_8);$

$\bullet$ the counit

$\varepsilon(x_1)=4$, $\varepsilon(x_6)=2$, $\varepsilon(x_2)=\varepsilon(x_3)=2$, $\varepsilon(x_7)=\varepsilon(x_8)=1$,

$\varepsilon(x_4)=\varepsilon(x_5)=\varepsilon(x_9)=\varepsilon(x_{10})=0$;

Then we define a map $\xi:H\rightarrow H$. Then it is an isomorphism of weak monoidal Hom-bialgebra if its matrix of the basis $\{x_{i}| i=1,2,\cdots 10\}$ takes the form
\[A=
\left(
  \begin{array}{cccccccccc}
    1 & 0 & 0 & 0 & 0 & 0 & 0 & 0 & 0 & 0 \\
    0 & 1 & 0 & 0 & 0 & 0 & 0 & 0 & 0 & 0 \\
    0 & 0 & 1 & 0 & 0 & 0 & 0 & 0 & 0 & 0 \\
    0 & 0 & 0 & \lambda & 0 & 0 & 0 & 0 & 0 & 0 \\
    0 & 0 & 0 & 0 & \lambda & 0 & 0 & 0 & 0 & 0 \\
    0 & 0 & 0 & 0 & 0 & 1 & 0 & 0 & 0 & 0 \\
    0 & 0 & 0 & 0 & 0 & 0 & 1 & 0 & 0 & 0 \\
    0 & 0 & 1 & 0 & 0 & 0 & 0 & 1 & 0 & 0 \\
    0 & 0 & 0 & 0 & 0 & 0 & 0 & 0 & \lambda & 0 \\
    0 & 0 & 0 & 0 & 0 & 0 & 0 & 0 & 0 & \lambda \\
  \end{array}
\right)
\]
Then it is a direct computation to check that $(H,\xi)$ is a weak monoidal Hom-bialgebra.\\

{\bf Definition 2.11.}  A weak monoidal Hom-bialgebra $(H,\xi)$ is called a {\sl weak monoidal Hom-Hopf algebra} if $H$ endowed with a $k$-linear map $S$ (the \textsl{antipode}), such that for any $h,g \in H$, the following conditions hold:

$
\left.\begin{array}{l}
(1)~~S\circ \xi = \xi\circ S;\\
(2)~~h_1S(h_2) = \varepsilon_t(h),~~S(h_1)h_2 = \varepsilon_s(h);\\
(3)~~S(hg) = S(g)S(h),~~S(1_H) = 1_H;\\
(4)~~\Delta(S(h)) = S(h_2) \o S(h_1),~~\varepsilon\circ S = \varepsilon.
\end{array}\right.
$
\\

{\bf Proposition 2.12.}  $H$ is a weak monoidal Hom-Hopf algebra, then for any $h \in H$, the following equalities hold:
\begin{eqnarray}
&&\mbox{(i)}~~~\varepsilon_s(h_1)S(h_2) = \xi^{-1}(S(h)),~~S(h_1)\varepsilon_t(h_2) = \xi^{-1}(S(h));\\
&&\mbox{(ii)}~~\varepsilon_t(h) = S(\widehat{\varepsilon_s}(h)),~~\varepsilon_s(h) = S(\widehat{\varepsilon_t}(h));\\
&&\mbox{(iii)}~\varepsilon_t(h_1) \o h_2 = S(1_1) \o 1_2\xi^{-2}(h),~~h_1 \o \varepsilon_s(h_2) = \xi^{-2}(h)1_1 \o S(1_2).
\end{eqnarray}
{\bf Proof. } We only check the first identities of each one. \\
(i).
\begin{eqnarray*}
\varepsilon_{s}(h_{1})S(h_{2})&=& 1_1\varepsilon(h_11'_11_2)S(\xi^{-1}(h_2)1'_2)\\
&=&1_1S(\xi^{-2}(h)1_2)=\xi^{-2}(S(h))1=\xi^{-1}(S(h));
\end{eqnarray*}
(ii)
\begin{eqnarray*}
 \varepsilon_t(h)&=& \varepsilon(1_1\varepsilon_{t}(h))1_2\stackrel{(2.21)}{=}\varepsilon(\varepsilon_{t}(h)1_1)1_2\\
 &\stackrel{(2.15)}{=}&\varepsilon(\varepsilon_{s}(h_{1})S(h_{2})1_1)1_2= \varepsilon(S(h)1_1)1_2\\
 &=& \varepsilon(S(h)S(1_2))S(1_1)= S(1_1)\varepsilon(1_2h)= S(\widehat{\varepsilon_s}(h));
\end{eqnarray*}
(iii) Easy to get from (ii).
$\hfill \blacksquare$
\\

Note that, a monoidal Hom-bialgebra is Hom-bialgebra if and only if the Hom-structure map $\xi$ satisfies
$\xi \circ \xi = id$. Moreover, through a direct computation, we can get that there is a one to one correspondence between the collection the monoidal Hom-bialgebras over a commutative ring $k$, and the collection of the unital Hom-bialgebra over $k$ which Hom-structure map is a bijection. \\

Recall from \cite{ZW}, a weak Hom-bialgebra is a Hom-algebra and a Hom-coalgebra with the compatible conditions as follows:

\[\left.\begin{array}{l}
(1)~~\Delta(ab)  = \Delta(a)\Delta(b);\\
(2)~~\varepsilon((ab)c) = \varepsilon(ab_1)\varepsilon(b_2c),~~\varepsilon(a(bc)) = \varepsilon(ab_2)\varepsilon(b_1c);\\
(3)~~(\Delta \o id_H)\Delta(1_H) = 1_1 \o 1_21'_{1} \o 1'_{2},~~(id_H \o \Delta)\Delta(1_H) = 1_1 \o1'_{1}1_2 \o 1'_{2}.
\end{array}\right.
\]
\\

More precisely, we can obtain the following relationship between the collection of weak monoidal Hom-Hopf algebras and weak Hom-Hopf algebras.\\

{\bf Proposition 2.13.}
If $(H,\xi,m,1_H,\Delta,\varepsilon,S)$ is a weak monoidal Hom-Hopf algebra, then ${}^\xi H=(H,\xi,m,1_H,\Delta\circ\xi^2,\varepsilon,S)$ is a weak Hom-Hopf algebra. Conversely, if $(B$, $\xi_B$, $m$, $1_B, \Delta, \varepsilon, S)$ is a weak Hom-Hopf algebra and $\xi_B$ is invertible, then ${}_\xi B=(B,\xi_B,m,1_B,\Delta\circ\xi_B^{-2},\varepsilon,S)$ is a weak monoidal Hom-Hopf algebra.

{\bf Proof.} Firstly we denote $\Delta\circ\xi^2(h)$ by $h_{[1]} \o h_{[2]}$. Then we have
\begin{eqnarray*}
\xi(h_{[1]} )\o h_{[2][1]} \o h_{[2][2]}&=& \xi(h_{[1]}) \o \xi^2( h_{[2]1})\o \xi^2( h_{[2]2})\\
&=& \xi^3(h_{1}) \o \xi^4( h_{21})\o \xi^4( h_{22})\\
&=& \xi^4(h_{11}) \o \xi^4( h_{12})\o \xi^3( h_{2})\\
&=&h_{[1][1]}\o h_{[1][2]} \o \xi(h_{[2]} ),
\end{eqnarray*}
and
\begin{eqnarray*}
\varepsilon(h_{[1]})h_{[2]}&=&\varepsilon(h_{1})\xi^2( h_{2})\\
&=&\xi(h)\\
&=&h_{[1]}\varepsilon(h_{[2]}),
\end{eqnarray*}
which implies $({}^\xi H,\xi)$ is a Hom-coalgebra. Obviously $({}^\xi H,\xi)$ is a Hom-algebra.

Secondly, we have
\begin{eqnarray*}
\varepsilon((x y)z)=\varepsilon(x y_1)\varepsilon(y_2z)=\varepsilon(x y_{[1]})\varepsilon(y_{[2]}z),
\end{eqnarray*}
and the other conditions are easily to get. The proof of the opposite statement is left to readers.

$\hfill \blacksquare$\\

{\bf Proposition 2.14.} We find that (\cite{ZW}, Proposition 2.9) should be like:

For any given weak bialgebra $(H,\mu,1_H,\Delta,\varepsilon)$, suppose that $\a:H\rightarrow H$ is both a morphism of algebras preserving unit and a morphism of coalgebras preserving counit. Thus we can define a new multiplication $\overline{\mu}:=\a\circ\mu$, and a new comultiplication $\overline{\Delta}: =\Delta\circ \a$,
then $H^\a=(H,\a,\overline{\mu},\eta,\overline{\Delta},\varepsilon)$ is a weak Hom-bialgebra if and only if $\a$ satisfies $\a(1_1)\o 1_2 = \Delta(1_H)$.

{\bf Proof.} $\Leftarrow$: We denote $\overline{\mu}(a \o b)$ by $a \ci b$, $\overline{\Delta}(c) = c_{[1]} \o c_{[2]}$ for any $a,b,c \in H$.

Firstly, we have
$$\aligned
(a \ci b) \ci \a(c) &= \a(\a(a) \a(b))\a^2 (c)\\
&=\a(\a(a)\a(b c)) = \a(a)\ci(b\ci c).
\endaligned$$
Since $\a(1_H) = 1_H$, thus $(H,\a,\overline{\mu},1_H)$ is a Hom-algebra.

Similarly, $(H,\a,\overline{\Delta},\varepsilon)$ is a Hom-coalgebra.

Secondly, from $\a(1_1)\o 1_2 = \Delta(1_H)$, we have $1_1\o \a(1_2) = \Delta(1_H)$. Thus we can get
$1_{[1]} \o 1_{[2][1]} \o 1_{[2][2]} = 1_{[1]} \o 1'_{[1]}1_{[2]} \o 1'_{[2]}$ and $1_{[1][1]} \o 1_{[1][2]} \o 1_{[2]} = 1_{[1]} \o1_{[2]}1'_{[1]} \o 1'_{[2]}$.

Thirdly, since $\a(1_1)\o 1_2 = \Delta(1_H)$, we immediately get
$$\aligned
\varepsilon(\a(a)b)&= \varepsilon(\a(a)1_1)\varepsilon(1_2b)\\
&=\varepsilon(\a(a)\a(1_1))\varepsilon(1_2b)  = \varepsilon(a1_1)\varepsilon(1_2b)\\
&=\varepsilon(ab),
\endaligned$$
then we have
$$\aligned
\varepsilon((a\ci b) \ci c)&= \varepsilon(\a(\a(ab)c)) = \varepsilon(\a^2(a)\a^2(b_1))\varepsilon(\a^2(b_2)\a(c))\\
&=\varepsilon(\a(a)\a^2(b_1))\varepsilon(\a(\a(b_2)_2c)) \\
&= \varepsilon(a \ci b_{[1]})\varepsilon(b_{[2]}\ci c),
\endaligned$$
and similarly we can obtain $\varepsilon(a\a(b))=\varepsilon(ab)$ through $1_1\o \a(1_2) = \Delta(1_H)$. Thus we have $\varepsilon(a\ci (b \ci c)) = \varepsilon(a \ci b_{[2]})\varepsilon(b_{[1]}\ci c)$.

Finally, we check that
$$\aligned
(a\ci b)_{[1]} \o (a\ci b)_{[2]}&=\a({\a(a)}_1{\a(b)}_1) \o \a({\a(a)}_2{\a(b)}_2) \\
&=a_{[1]}\ci b_{[1]} \o a_{[2]}\ci b_{[2]},
\endaligned$$
which means that $H^\a=(H,\a,\overline{\mu},1_H,\Delta,\varepsilon)$ is a weak Hom-bialgebra.

$\Rightarrow$: Straightforward.
$\hfill \blacksquare$
 \\

{\bf Remark.} Based on the above proposition, it is easy to know that every 2-dimensional weak Hom-bialgebras have trivial structures. That means if $H_2$ is a $\Bbbk$-space with a basis ${I,E}$, and the following structures

$\bullet $ the multiplication
$$\begin{array}{|c|c|c|}
\hline H_2 & I  & E  \\
\hline I   & I  & E \\
  E  & E  & E   \\
 \hline
 \end{array}$$

$\bullet $ the comultiplication
$$\Delta(I) = (I-E) \o(I-E)+E \o E,~~\Delta(E) = E \o E;$$

$\bullet $ the counit
$$\varepsilon(I) = 2,~~\varepsilon(E) = 1.$$

Then the automorphism of $H_2$ is identity map. This means (\cite{ZW}, Example 2.12) is a trivial weak bialgebra.
\\

\section*{3. Weak $(\alpha, \beta)$- Yetter-Drinfeld monoidal Hom-modules}
\def\theequation{3. \arabic{equation}}
\setcounter{equation} {0} \hskip\parindent

 In this section, we will define the notion of a Yetter-Drinfeld module over
 a weak monoidal Hom-Hopf algebra that is twisted by two weak monoidal
 Hom-Hopf algebra automorphisms as well as the notion of a weak monoidal Hom-entwining structure
 and how to obtain such structure from automorphisms of weak monoidal Hom-Hopf algebras.
\\

In what follows, let  $(H,\xi)$ be a weak monoidal Hom-Hopf algebra with
the bijective antipode $S$  and let ${\sl Aut}_{wmHH}(H)$ denote the set
  of all automorphisms of a weak monoidal Hopf algebra $H$.
\\

{\bf Definition 3.1.} Let $\alpha, \beta\in  {\sl Aut}_{wmHH}(H)$.
A weak left-right {\sl $(\alpha, \beta)$-Yetter-Drinfeld Hom-module} over $(H,\xi)$
is a vector space $M$ such that:
\begin{eqnarray*}
 &(1)&(M,\cdot,\mu) \mbox{ is a left }H \mbox{-Hom-module;}\\
 &(2)& (M,\rho,\mu) \mbox{ is a right }H \mbox{-Hom-comodule;}\\
 &(3)& \rho \mbox{ and }\cdot \mbox{ satisfy
  the following compatibility condition: }
   \end{eqnarray*}
 \begin{equation}
\r (h\c m)=\xi(h_{21})\c m_{(0)}\o (\beta(h_{22})\xi^{-1}(m_{(1)}))\alpha(S^{-1}(h_{1})),
\end{equation}
for all $h\in H$ and $m\in M$. We denote by
$ _{H}\mathcal{WMHYD}^{H}(\alpha, \beta)$
the category of weak left-right $(\alpha, \beta)$-Yetter-Drinfeld
Hom-modules, morphisms
being $H$-linear $H$-colinear maps.
\\

 {\bf Remark.} Note that, $\alpha$ and $\beta$ are bijective,
 Hom-algebra morphisms,
 Hom-coalgebra morphisms,
 and commute with $S$ and $\xi$.
 \\

{\bf Proposition 3.2.} One has that Eq.(3.1) is equivalent to
 the following equations:
\begin{eqnarray}
\rho(m)=m_{(0)}\otimes m_{(1)}\in M\otimes_{t_{\beta}} H \stackrel{\triangle}{=} (1_{1}\otimes \beta(1_{2}))\cdot (M\otimes H), \forall m\in M,\\
h_{1}\c m_{(0)}\o \beta(h_{2})m_{(1)}=\mu((h_{2}\c\mu^{-1} (m))_{(0)}) \o (h_{2}\c
 \mu^{-1}(m))_{(1)}\alpha(h_{1}).
\end{eqnarray}

{\bf Proof.} Eq.(3.1)$\Longrightarrow$ Eq.(3.2, 3.3). We first note that
\begin{eqnarray*}
  m_{(0)}\otimes m_{(1)}
  &=& \xi(1_{21})\cdot \mu^{-1}(m)_{(0)}\otimes (\beta(1_{22})\xi^{-2}(m_{(1)}))\alpha(S^{-1}(1_{1}))\\
  &=& 1_{21}\cdot \mu^{-1}(m)_{(0)}\otimes ((\beta(1_{22})\xi^{-2}(m_{(1)}))\alpha(S^{-1}(1_{1})))\\
  &=& 1'_{1}\cdot (1_{2}\cdot \mu^{-2}(m)_{(0)})\otimes \beta(1'_{(2)})(\xi^{-2}(m_{(1)})\alpha(S^{-1}(1_{1})))\in M\otimes_{t_{\beta}} H.
\end{eqnarray*}

Then we do calculation as follows:
\begin{eqnarray*}
 &&\mu((h_{2}\c \mu^{-1}(m))_{(0)}) \o (h_{2}\c\mu^{-1}( m))_{(1)}\alpha(h_{1})\\
&\stackrel{(3.1)}{=}&  \mu(\xi(h_{221})\cdot \mu^{-1}(m)_{(0)})\otimes(( \beta(h_{222})\xi^{-2}(m_{(1)}))\alpha(S^{-1}(h_{21})))\alpha(h_{1})\\
 &=& \xi^{-1}(h_{1}1_{2})\cdot m_{(0)} \otimes \xi^{-1}(\beta(h_{2})m_{1})\alpha(S^{-1}(1_{1}))\\
 &=& h_{1}\cdot (1_{2}\cdot \mu^{-1}(m_{0}))\otimes \xi^{-1}(\beta(h_{2})(\beta(1_{3})\xi^{-1}(m_{1})))\alpha(S^{-1}(1_{1}))\\
 &=& h_{1}\cdot m_{(0)}\otimes \beta(h_{2})m_{(1)}.
\end{eqnarray*}

 For Eq.(3.2, 3.3) $\Longrightarrow$ Eq.(3.1), we have
\begin{eqnarray*}
&& \xi(h_{21})\cdot m_{(0)}\otimes_{t_{\beta}} (\beta(h_{22})\xi^{-1}(m_{(1)}))\alpha(S^{-1}(h_{1}))\\
&\stackrel{(3.3)}{=}& \mu((\xi(h_{22})\cdot \mu^{-1}(m))_{(0)})\otimes_{t_{\beta}} \xi^{-1}((\xi(h_{22})\cdot \mu^{-1}(m))_{(1)}\alpha(\xi(h_{21})))\alpha(S^{-1}(h_{1}))\\
&=& \mu(((1_{2}\xi^{-2}(h))\cdot \mu^{-1}(m))_{(0)})\otimes_{t_{\beta}} ((1_{2}\xi^{-2}(h))\cdot \mu^{-1}(m))_{(1)}\alpha(1_{1})\\
&=& \mu((1_{2}\cdot \mu^{-1}(\mu^{-1}(h\cdot m)))_{(0)})\otimes_{t_{\beta}} (1_{2}\cdot \mu^{-1}(\mu^{-1}(h\cdot m)))_{(1)}\alpha(1_{1})\\
&\stackrel{(3.3)}{=}& 1_{1}\cdot \mu^{-1}(h\cdot m)_{(0)}\otimes_{t_{\beta}} \beta(1_{2})\xi^{-1}(h\cdot m)_{(1)}\\
&\stackrel{(3.2)}{=}& 1'_{1}\cdot(1_{1}\cdot \mu^{-1}(h\cdot m)_{(0)})\otimes \beta(1'_{2})\cdot(\beta(1_{2})\xi^{-1}(h\cdot m)_{(1)})\\
&=& (h\cdot m)_{(0)}\otimes (h\cdot m)_{(1)}.
\end{eqnarray*}
This finishes the proof.
\hfill $\blacksquare$
\\

{\bf Definition 3.3.} A weak left-right {\sl monoidal Hom-entwining structure} is a triple
 $(H, \, C,\, \psi)$, where $(H, \xi)$ is a monoidal Hom-algebra
 and $(C,\gamma)$ is a monoidal Hom-coalgebra with a linear map
 $\psi: H\otimes C\rightarrow H\otimes C,
 \,\,\,\,h\otimes c\mapsto _{\psi}h\otimes c^{\psi}$
 satisfying the following conditions:
\begin{equation}
 _{\psi}(hg)\otimes c^{\psi}= _{\phi}h_{\psi}g\otimes \gamma(\gamma^{-1}(c)^{\psi\phi}),
\end{equation}
\begin{equation}
 _{\psi}1_{H}\otimes c^{\psi}= \varepsilon(c_{1}^{\psi}) {_{\psi}1_{H}}\otimes \gamma(c_{2}),
\end{equation}
\begin{equation}
 _{\psi}h\otimes \Delta(c^{\psi})= \xi(_{\phi\psi}\xi^{-1}(h))\otimes (c^{\psi}_{1}\otimes c^{\phi}_{2}),
\end{equation}
\begin{equation}
 \varepsilon(c^{\psi}) _{\psi}h= \varepsilon(\gamma^{-1}(c)^{\psi})h(_{\psi}1_{H}),
\end{equation}\\

 Over a weak monoidal Hom-entwining structure $(H, \, C,\, \psi)$,
 a left-right weak monoidal entwined Hom-module $M$ is both a right
 $C$-Hom-comodule and a left $H$-Hom-module such that
 $$ \rho^{M}(h\cdot m)={_{\psi}\xi^{-1}(h)}\cdot m_{(0)}\otimes \gamma(m_{(1)})^{\psi}$$
 for all $h\in H$ and $m\in M$. We denote the category of all
 monoidal entwined Hom-modules over $(H,\,C,\,\psi)$
 by $_{H}\mathcal{M}^{C}(\psi)$.
 \\

 Let $(H,\xi)$ be a weak monoidal Hom-Hopf algebra with $S$,
 and define a linear map
 $$\psi(\alpha, \beta):H\otimes H\rightarrow H\otimes H,\
 \ a\otimes c\mapsto {_{\psi}a}\otimes c^{\psi}
 =\xi^{2}(a_{21})\otimes (\beta(a_{22})\xi^{-2}(c))\alpha(S^{-1}(a_1)),$$
 for all $\alpha, \beta\in  {\sl Aut}_{wmHH}(H)$.
\\

{\bf Proposition 3.4.} With notations above,
 $(H, H,\psi(\alpha, \beta))$ is a weak monoidal Hom-entwining structure
 for all $\alpha, \beta\in  {\sl Aut}_{wmHH}(H)$.

 {\bf Proof.} We need to prove that Eqs.(3.4-3.7) hold.
 First, it is straightforward to check Eqs.(3.5) and (3.7).
 In what follows, we only verify  Eqs.(3.4) and (3.6).
 In fact, for all $a,b,c\in H$, we have
 \begin{eqnarray*}
&&_{\phi}a_{\psi}b\otimes \xi(\xi^{-1}(c)^{\psi\phi})\\
&=&\xi^2(a_{21})_{\psi}b \otimes
\xi((\beta(a_{22})\xi^{-2}(\xi^{-1}(c)^{\psi}))\alpha S^{-1}(a_1))\\
&=&\xi^2(a_{21}b_{21})
\otimes(\beta(a_{22}b_{22})\xi^{-2}(c))(\alpha S^{-1}(b_1)\alpha S^{-1}(a_1))
= _{\psi}(ab)\otimes c^{\psi},
\end{eqnarray*}
and Eq.(3.4) is proven.

For all $a\in H$, we have
\begin{equation}
a_{1}\otimes a_{211}\otimes a_{2121}\otimes a_{2122}\otimes a_{22}
=\xi(a_{11})\otimes \xi^{-1}(a_{12})\otimes \xi^{-2}(a_{21})
\otimes \xi^{-1}(a_{221})\otimes \xi(a_{222})
\end{equation}

As for Eq.(3.6), we compute:
\begin{eqnarray*}
&& \xi(_{\phi\psi}\xi^{-1}(a))\otimes (c^{\phi}_{1}\otimes c^{\psi}_{2})\\
&=&\xi(\xi^{2}((_{\psi}\xi^{-1}(a))_{21}))\otimes
((\beta((_{\psi}\xi^{-1}(a))_{22})\xi^{-2}(c_{1}))\alpha S^{-1}((_{\psi}\xi^{-1}(a))_{1})
\otimes c^{\psi}_{2})\\
&=&\xi^{2}(a_{21})\otimes
((\beta(a_{221})\xi^{-2}(c_{1}))\alpha S^{-1}(a_{12})
\otimes (\beta(a_{222})\xi^{-2}(c_{2}))\alpha S^{-1}(a_{11}))\\
&=&_{\psi}a\otimes \Delta(c^{\psi}).
\end{eqnarray*}
and Eq.(3.6) is proven.

This finishes the proof. \hfill $\blacksquare$
\\

{\bf Remark. } By Proposition above, we have a weak monoidal entwined
  Hom-module category $ _{H}\mathcal{M}^{H}(\psi(\alpha, \beta))$
 over $(H, H,\psi(\alpha, \beta))$ with $\alpha, \beta\in{\sl Aut}_{wmHH}(H)$.
 In this case, for all $M\in _{H}\mathcal{M}^{H}(\psi(\alpha, \beta))$,
 we have
  $$\r (h\c m)=\xi(h_{21})\c m_{(0)}
  \o (\beta(h_{22})\xi^{-1}(m_{(1)}))\alpha(S^{-1}(h_{1})),$$
  for all $h\in H,m\in M$. Thus means that
  $_{H}\mathcal{M}^{H}(\psi(\alpha, \beta))=
   _{H}\mathcal{WMHYD}^{H}(\alpha, \beta)$
  as categories.
  \\

\section*{4. A BRAIDED $T$-CATEGORY $\mathcal {WMHYD}(H)$}
\def\theequation{4. \arabic{equation}}
\setcounter{equation} {0} \hskip\parindent

In this section, we will construct a class of new braided $T$-categories
 $\mathcal {WMHYD}(H)$ over any weak monoidal Hom-Hopf algebra $(H,\xi)$ with bijective antipode.\\

 Let $(M,\mu) \in {_{H}}\mathcal {WMHYD}^{H}(\alpha,\beta)$
 ,$(N,\nu)\in {_{H}}\mathcal {WMHYD}^{H}(\gamma,\delta)$, with
 $\alpha,\beta,\gamma,\delta \in {\sl Aut}_{wmHH}(H)$. Define $M\otimes_{t_{\gamma^{-1}\beta}} N= (1_{1}\otimes \gamma^{-1}\beta(1_{2}))\cdot (M\otimes N)$.
\\

{\bf Proposition 4.1.} If $(M,\mu) \in {_{H}}\mathcal {WMHYD}^{H}(\alpha,\beta)$
 and $(N,\nu)\in {_{H}}\mathcal {WMHYD}^{H}(\gamma,\delta)$, with
 $\alpha,\beta,\gamma,\delta \in {\sl Aut}_{wmHH}(H)$, then $(M \otimes_{t_{\gamma^{-1}\beta}} N,\mu\otimes \nu)
 \in {_{H}}\mathcal {WMHYD}^{H}(\alpha\gamma, \delta\gamma^{-1}\beta\gamma)$ with structures as follows:
\begin{eqnarray*}
&&h\c (m \otimes n)=\gamma (h_{1})\c m \otimes \gamma^{-1}\beta\gamma(h_{2})\c n,\\
&&m\otimes n \mapsto (m_{(0)}\otimes n_{(0)})\otimes n_{(1)}m_{(1)}.
\end{eqnarray*}
for all $m\in M,n\in N$ and $h\in H.$
\\

{\bf Proof.} Let $h,g\in H$ and $m\otimes n\in M\otimes_{t_{\gamma^{-1}\beta}} N$. We can prove
$(hg)\cdot (m\otimes n)=\xi(h)\cdot (g\cdot (\mu^{-1}(m)\otimes \nu^{-1}(n))$ straightforwardly, and
\begin{eqnarray*}
 1\cdot (m\otimes_{t_{\gamma^{-1}\beta}} n)&=& \gamma(1'_{1})\cdot (1_{1}\cdot m)\otimes \gamma^{-1}\beta\gamma(1'_{2})\cdot (\gamma^{-1}\beta(1_{2})\cdot n)\\
 &=& (\gamma(1'_{1})1_{1})\cdot \mu(m)\otimes (\gamma^{-1}\beta\gamma(1'_{2})\gamma^{-1}\beta(1_{2}))\cdot \nu(n)\\
 &=& 1_{1}\cdot \mu(m)\otimes \gamma^{-1}\beta(1_{2})\cdot \nu(n)\\
 &=& \mu(m)\otimes_{t_{\gamma^{-1}\beta}} \nu(n).
\end{eqnarray*}
This shows that $(M \otimes_{t_{\gamma^{-1}\beta}} N,\mu\otimes \nu)$ is a left $H$-module
and the right $H$-comodule condition is straightforward to check.

Next, we compute the compatibility condition as follows:
\begin{eqnarray*}
&&(h\c (m\otimes n))_{(0)}\otimes (h\c (m\otimes n))_{(1)}\\
 &=&((\gamma(h_{1})\c m)_{(0)} \otimes (\gamma^{-1}\beta\gamma(h_{2})\c n)_{(0)})\otimes
 (\gamma^{-1}\beta\gamma(h_{2})\c n)_{(1)} (\gamma(h_{1})\c m)_{(1)}\\
&\stackrel{(3.1)}{=}&(\gamma \xi (h_{121})\c m_{(0)} \otimes  \gamma^{-1}\beta\gamma\xi(h_{221})\c n_{(0)})\otimes
((\delta\gamma^{-1}\beta\gamma(h_{222})\xi^{-1}(n_{(1)}))\\
&& \gamma S^{-1}\gamma^{-1}\beta\gamma(h_{21}))((\beta\gamma(h_{122})\xi^{-1}(m_{(1)}))\alpha S^{-1}\gamma(h_{11}))\\
&=&(\gamma (h_{12})\c m_{(0)} \otimes  \gamma^{-1}\beta\gamma\xi(h_{221})\c n_{(0)})\otimes
(\delta\gamma^{-1}\beta\gamma\xi(h_{222})n_{(1)}) \\
&&((\beta\gamma\xi^{-1}(\varepsilon(h_{21}1_{1})1_{2})\xi^{-1}(m_{(1)}))S^{-1}\alpha\gamma(h_{11}))\\
&=&(\gamma (h_{12})\varepsilon(h_{21}1_{1})\c m_{(0)} \otimes  \gamma^{-1}\beta\gamma\xi(h_{221})\c n_{(0)})\otimes
(\delta\gamma^{-1}\beta\gamma\xi(h_{222})n_{(1)}) \\
&&((\beta\gamma\xi^{-1}(1_{2})\xi^{-1}(m_{(1)}))S^{-1}\alpha\gamma(h_{11}))\\
&=& (\gamma\xi^{-2}(h_{12}1_{2}))\cdot m_{(0)}\otimes \gamma^{-1}\beta\gamma(h_{21})\cdot n_{(0)})
\otimes (\delta\gamma^{-1}\beta\gamma(h_{22})n_{(1)})\\
&& ((\beta\gamma(1_{3})\xi^{-1}(m_{(1)}))S^{-1}\alpha\gamma\xi^{-1}(h_{11}1_{1}))\\
&=& \gamma(h_{12})\cdot(1_{2}\cdot \mu^{-1}(m_{(0)}))\otimes \gamma^{-1}\beta\gamma(h_{21})\cdot n_{(0)}\otimes (\delta\gamma^{-1}\beta\gamma(h_{22})n_{(1)})\\
&&(((\beta(1_{3})\xi^{-2}(m_{(1)}))\alpha S^{-1}(1_{1}))\alpha\gamma S^{-1}(h_{11}))\\
&=&(\gamma \xi(h_{211})\c m_{(0)} \otimes  \gamma^{-1}\beta\gamma\xi(h_{212})\c n_{(0)})\otimes
(\delta\gamma^{-1}\beta\gamma(h_{22})\xi^{-1}(n_{(1)}m_{(1)}))\\
&&S^{-1}\alpha\gamma(h_{1})\\
&=&\xi(h_{21})\c (m\otimes n)_{(0)}\otimes \delta\gamma^{-1}\beta\gamma(h_{22})\xi^{-1}(m\otimes
n)_{(1)}\alpha\gamma(S^{-1}(h_{1})).
\end{eqnarray*}
Thus $(M \otimes_{t_{\gamma^{-1}\beta}} N,\mu\otimes \nu)
 \in {_{H}}\mathcal {WMHYD}^{H}(\alpha\gamma, \delta\gamma^{-1}\beta\gamma)$.
 \hfill $\blacksquare$
\\

{\bf Remark. } Note that, if $(M,\mu)\in {_{H}}\mathcal
{WMHYD}^{H}(\alpha,\beta),\,\, (N,\nu) \in {_{H}}\mathcal {WMHYD}^{H}(\gamma,\delta)$ \\
and
$(P,\varsigma) \in {_{H}}\mathcal {WMHYD}^{H}(s, t)$, then the associativity constraint $a_{M,N,P}$ is
\begin{eqnarray*}
a_{M,N,P}:(M\otimes_{t_{\gamma^{-1}\beta}} N)\otimes_{t_{s^{-1}\delta\gamma^{-1}\beta\gamma}} P&\rightarrow&
M \otimes_{t_{s^{-1}\gamma^{-1}\beta}}(N\otimes_{t_{s^{-1}\delta}} P)\\
(m\otimes n)\otimes p &\mapsto&
\mu(m) \otimes(n\otimes\varsigma^{-1}(p))
\end{eqnarray*}
where $(M\otimes_{t_{\gamma^{-1}\beta}} N)\otimes_{t_{s^{-1}\delta\gamma^{-1}\beta\gamma}} P\in {_{H}\mathcal {WMHYD}}^{H}(\alpha\gamma s,
ts^{-1}\delta\gamma^{-1}\beta\gamma s)$.
\\

 Denote $G={\sl Aut}_{wmHH}(H)
 \times {\sl Aut}_{wmHH}(H)$
 a group with multiplication as follows:
  for all $\alpha,\beta, \gamma, \delta \in {\sl Aut}_{wmHH}(H)$,
\begin{equation}
(\alpha,\beta)\ast (\gamma, \delta)=(\alpha\gamma, \delta\gamma^{-1}\beta\gamma).
\end{equation}
 The unit of this group is $(id,id)$ and $(\alpha,\beta)^{-1}=(\alpha^{-1},
 \alpha\beta^{-1}\alpha^{-1})$.\\

The above proposition means that if
$M \in {_{H}}\mathcal {WMHYD}^{H}(\alpha,\beta)$
 and $N\in {_{H}}\mathcal {WMHYD}^{H}(\gamma, \delta)$, then
 $M \otimes N \in {_{H}}\mathcal {WMHYD}^{H}((\alpha,\beta)\ast (\gamma, \delta)).$
\\

{\bf Proposition 4.2.} The associativity constraints of monoidal category ${_{H}}\mathcal{WMHYD}^{H}$ are described as above.
The left and right unit constraints $l_{N}:H_{t}\otimes_{t_{\gamma^{-1}}} N\rightarrow N$
and $r_{M}:M\otimes_{t_{\beta}} H_{t}\rightarrow M$ with $(N,\nu)\in {_{H}}\mathcal {WMHYD}^{H}(\gamma,\delta)$ and $(M,\mu)\in {_{H}}\mathcal {WMHYD}^{H}(\alpha,\beta)$ and their
inverses are given by the formulas
\begin{eqnarray*}
 &&l_{N}(x\otimes_{t_{\gamma^{-1}}}n)=\gamma^{-1}(x)\cdot n,\quad l^{-1}_{N}(n)= \varepsilon_{t}(1_{1})\otimes_{t_{\gamma^{-1}}} \gamma^{-1}(1_{2})\cdot \nu^{-2}(n);\\
 &&r_{M}(m\otimes_{t_{\beta}}x)= \widehat{\varepsilon_{s}}(\beta^{-1}(x))\cdot m,\quad r^{-1}_{M}(m)=\beta^{-1}(1_{1})\cdot \mu^{-2}(m) \otimes_{t_{\beta}} 1_{2}.
\end{eqnarray*}
for all $ x\in H_{t}$, $n\in N$ and $m\in M$.

{\bf Proof.}
Observe that $H_{t}\in {_{H}}\mathcal {WMHYD}^{H}(id,id)$, with left $H$-action $h\cdot x=\varepsilon_{t}(hx)$ and right $H$-coaction $\rho(x)=1_{2}\otimes S^{-1}(x1_{1})$, for all $h\in H$, here $\rho$ is the comodule structure map.

It is easy to get that $H_{t}$ is a left $H$-module and a right $H$-comodule under"$\cdot$" and"$\rho$". We just check $H_{t}\in {_{H}}\mathcal {WMHYD}^{H}(id,id)$.

On the one hand,
\begin{eqnarray*}
\rho(h\cdot x) &=& \rho(\varepsilon_{t}(hx))\\
&=& 1_{2}\otimes S^{-1}(1_{1})S^{-1}((\xi^{-2}(h_{1})x)S(\xi^{-1}(h_{2})))\\
&=& 1_{2}\otimes S^{-1}(1_{1})(\xi^{-1}(h_{2})S^{-1}(\xi^{-2}(h_{1})x)).
\end{eqnarray*}

On the other hand,
\begin{eqnarray*}
\rho(h\cdot x)&=& \xi^(h_{21})\cdot 1_{2}\otimes (h_{22}\xi^{-1}(S^{-1}(x1_{1})))S^{-1}(h_{1})\\
&=& \varepsilon_{t}(\xi(h_{21})1_{2})\otimes (h_{22}(S^{-1}(x1_{1})))S^{-1}(h_{1})\\
&=& 1_{2}\otimes \varepsilon(1_{1}h_{21})\xi(h_{22})S^{-1}(\xi^{-1}(h_{1})x)\\
&=& 1_{2}\otimes (\widehat{\varepsilon_{t}}(1_{1})\xi^{-1}(h_{2}))S^{-1}(\xi^{-1}(h_{1})x)\\
&=& 1_{2}\otimes S^{-1}(1_{1})(\xi^{-1}(h_{2})S^{-1}(\xi^{-2}(h_{1})x)).
\end{eqnarray*}

We just check the properties of $l_{N}$, the $r_{M}$ cases left to reader.
First we need to check $l_{N}$ is both $H$-linear and $H$-colinear.
Let $n\in N$. We have
\begin{eqnarray*}
l_{N}(h\cdot(x\otimes n))&=& l_{M}(\varepsilon_{t}(\alpha(h_{1})x)\otimes h_{2}\cdot n)\\
&=& \alpha^{-1}(\varepsilon_{t}(h_{1}x))\cdot(h_{2}\cdot n)\\
&=& (\xi^{-1}(h)\alpha^{-1}(x))\cdot \mu(n)\\
&=& h\cdot(\alpha^{-1}(x)\cdot n)\\
&=&h\cdot l_{M}(x\otimes n).
\end{eqnarray*}
So $l_{N}$ is $H$-linear, similar to the colinear case.

Next we check that
\begin{eqnarray*}
l_{N}(l^{-1}_{N}(n)) &=& l_{N}(\varepsilon_{t}(1_{1})\otimes \gamma^{-1}(1_{2})\cdot \nu^{-2}(n))\\
&=& \gamma^{-1}(\varepsilon_{t}(1_{1})1_{2})\cdot \nu^{-1}(n)\\
&=& 1\cdot \nu^{-1}(n)= n.
\end{eqnarray*}

\begin{eqnarray*}
l^{-1}_{N}(l_{N}(x\otimes n)) &=& l^{-1}_{N}(\gamma^{-1}\cdot n)\\
&=& \varepsilon_{t}(1_{1})\otimes \gamma^{-1}(1_{2})\cdot \nu^{-2}(\gamma^{-1}(x)\cdot n)\\
&=& \varepsilon_{t}(1_{1}x)\otimes \gamma^{-1}(1_{2})\cdot \nu^{-1}(n)\\
&=& (1'_{1}1_{1})\cdot x\otimes \gamma^{-1}(1'_{2}1_{2})\cdot n\\
&=& x\otimes n.
\end{eqnarray*}

Finally we need to prove the following diagram commute.\\

\centerline{
\xymatrix
{
(M\otimes_{t_{\beta}} H_{t})\otimes_{t_{\gamma^{-1}\beta}} N \ar[r]^-{a_{M,H_{t}, N}}\ar[d]_-{r_{M}\otimes id}
 & M\otimes_{t_{\gamma^{-1}\beta}}(H_{t}\otimes_{t_{\gamma^{-1}}} N)\ar[dl]^-{id\otimes l_{N}}\\
 M\otimes_{t_{\gamma^{-1}\beta}} N
}
}

\begin{eqnarray*}
(r_{M}\otimes id)((m\otimes x) \otimes n) &=& 1_{1}\cdot (\beta^{-1}(\widehat{\varepsilon_{s}}(x))\cdot m)\otimes \gamma^{-1}\beta(1_{2})\cdot n\\
&=& (1_{1}\beta^{-1}(\widehat{\varepsilon_{s}}(x))\cdot \mu(m) \otimes \gamma^{-1}\beta(1_{2})\cdot n\\
&=& \widehat{\varepsilon_{s}}(\beta^{-1}(x))_{1}\cdot \mu(m)\otimes \gamma^{-1}\beta(\varepsilon_{t}(\widehat{\varepsilon_{s}}(\beta^{-1}(x))_{2}))\cdot n\\
&=& 1_{1}\cdot \mu(m) \otimes \gamma^{-1}\beta(1_{2}\beta^{-1}(x))\cdot n\\
&=& 1_{1}\cdot \mu(m) \otimes \gamma^{-1}\beta(1_{2})\cdot (\gamma^{-1}(x)\cdot \nu^{-1}(n))\\
&=& (id\otimes l_{N})\circ a_{M,H_{t},N}((m\otimes x)\otimes n).
\end{eqnarray*}
This ends the proof.
\hfill $\blacksquare$\\

 {\bf Proposition 4.3.} Let $(N,\nu) \in {_{H}}\mathcal {WMHYD}^{H}(\gamma, \delta)$
 and $(\alpha,\beta)\in G$. Define ${}^{(\alpha,\beta)}N=N$ as
 vector space, with structures: for all $n\in N$ and $h\in H.$
$$h\rhd n=\gamma^{-1}\beta\gamma\alpha^{-1}(h)\c n,$$
\begin{equation}
n\mapsto n_{<0>}\otimes n_{<1>}=n_{(0)}\otimes \alpha\beta^{-1}(n_{(1)}).
\end{equation}
Then
$${}^{(\alpha,\beta)}N \in {_{H}}\mathcal {WMHYD}^{H}((\alpha,\beta)\ast (\gamma,\delta)\ast (\alpha,\beta)^{-1}).$$

{\bf Proof.} Obviously, the equations above define
a module and a comodule action.
In what follows,
we show the compatibility condition:
\begin{eqnarray*}
&&(h\rhd n)_{<0>}\otimes(h\rhd n)_{<1>}\\
&=&(\gamma^{-1}\beta\gamma\alpha^{-1}(h)\c n)_{(0)}\otimes \alpha\beta^{-1}((\gamma^{-1}\beta\gamma\alpha^{-1}(h)\c n)_{(1)})\\
&=&\gamma^{-1}\beta\gamma\alpha^{-1}\xi(h_{21})\c n_{(0)}\otimes (\alpha\beta^{-1}\delta\gamma^{-1}\beta\gamma\alpha^{-1}(h_{22})\alpha\beta^{-1}\xi^{-1}(n_{(1)}))
\alpha\gamma\alpha^{-1}S^{-1}(h_{1}))\\
&=&\xi(h_{21})\rhd n_{<0>}\otimes (\alpha\beta^{-1}\delta\gamma^{-1}\beta\gamma\alpha^{-1}(h_{22})\xi^{-1}(n_{<1>}))
\alpha\gamma\alpha^{-1}S^{-1}(h_{1}))
\end{eqnarray*}
for all $n\in N$ and $h\in H,$
that is ${}^{(\alpha,\beta)}N \in {_{H}}\mathcal {WMHYD}^{H}(\alpha\gamma\alpha^{-1},\alpha\beta^{-1}\delta\gamma^{-1}\beta\gamma\alpha^{-1})$
\hfill $\blacksquare$
\\

{\bf Remark.} Let $(M,\mu) \in {_{H}}\mathcal {WMHYD}^{H}(\alpha, \beta),
 \ \ (N,\nu) \in {_{H}}\mathcal
{WMHYD}^{H}(\gamma,\delta),\,\, \mbox {and}\, (s,t)\in G$. Then by the
 above proposition, we have:
 $$
 {}^{(\alpha, \beta)\ast (s,t)}N={}^{(\alpha, \beta)}({}^{(s,t)}N),
$$
as objects in $_{H}\mathcal {WMHYD}^{H}(\alpha s\gamma s^{-1}\alpha^{-1},
\alpha\beta^{-1}st^{-1}\delta\gamma^{-1}ts^{-1}\beta s\gamma s^{-1}\alpha^{-1})$ and
$$
{}^{(s,t)}(M\otimes N)= {}^{(s,t)}M \otimes {}^{(s,t)}N,
$$
as objects in $_{H}\mathcal {WMHYD}^{H}(s\alpha \gamma s^{-1},
 st^{-1}\delta\gamma^{-1}\beta\alpha^{-1}t\alpha \gamma s^{-1})$.
\\

 {\bf Proposition 4.4.} Let $(M,\mu) \in {_{H}}\mathcal {WMHYD}^{H}(\alpha,\beta)$
 and $(N,\nu)\in {_{H}}\mathcal {WMHYD}^{H}(\gamma,\delta)$, take
  ${}^{M}N={}^{(\alpha,\beta)}N$ as explained in Subsection 1.2.
  Define a map $c_{M, N}: M \otimes N \rightarrow {}^{M}N \otimes M$ by
\begin{equation}
 c_{M,N}(m\otimes
 n)=\nu(n_{(0)})\otimes \beta^{-1}(n_{(1)})\c \mu^{-1}(m).
\end{equation}
for all $m\in M,n\in N.$
Then $c_{M, N}$ is both an $H$-module map and an $H$-comodule map,
 and satisfies the following formulae
 (for $(P,\varsigma) \in {_{H}}\mathcal {WMHYD}^{H}(s,t)$):
\begin{equation}
a^{-1}_{{}^{M\otimes N}P, M, N}\circ c _{M\otimes N, P}\circ a^{-1}_{M, N, P}=(c _{M, {}^NP}\otimes
 id_N)\circ a^{-1}_{M, {}^NP, N}\circ (id _M\otimes c _{N, P})  ,
\end{equation}
\begin{equation}
 a_{{}^MN, {}^MP, M}
 \circ  c _{M, N\otimes P}\circ a_{M, N, P}=(id _{{}^MN }\otimes c _{M, P})\circ a_{{}^MN, M, P}\ci
 (c _{M, N}\otimes id_P).
\end{equation}

Furthermore, if $(M,\mu) \in {_{H}}\mathcal {WMHYD}^{H}(\alpha,\beta)$ and
$(N,\nu) \in{_{H}}\mathcal {WMHYD}^{H}(\gamma,\delta)$, then\\
 $c_{{}^{(s,t)}M,{^{(s,t)}N}}=c_{M,N},$
 for all $(s,t)\in G$.
\smallskip
\\

{\bf Proof.} First, we prove that $c_{M,N}$ is an $H$-module map.
 Take  $h\cdot(m\otimes n)=\gamma(h_{1})\c m\otimes \gamma^{-1}\beta\gamma(h_{2})\c n$
 and $h\cdot(n\otimes m)= \gamma^{-1}\beta\gamma(h_{1})\c n\otimes \beta^{-1}\delta\gamma^{-1}\beta\gamma(h_{2})\c m$
 as explained in Proposition 4.1.
 \begin{eqnarray*}
&&c_{M,N}(h\cdot(m\otimes n))\\
&=&\nu((\gamma^{-1}\beta\gamma(h_{2})\c n)_{(0)})\otimes
\beta^{-1}((\gamma^{-1}\beta\gamma(h_{2})\c n)_{(1)})\c\mu^{-1}(\gamma(h_{1})\c m)\\
&=&\nu(\gamma^{-1}\beta\gamma\xi(h_{221})\c n_{(0)})\otimes
\beta^{-1}(\delta\gamma^{-1}\beta\gamma\xi(h_{222}) n_{(1)})\\
&&\c ((\gamma S^{-1}\xi^{-1}(h_{21})
\gamma \xi^{-2}(h_{1}))\c \mu^{-1}(m))\\
&=& \nu(\gamma^{-1}\beta\gamma(\xi^{-2}(h_{1})1_{2})\cdot n_{(0)})\otimes \beta^{-1}(\delta\gamma^{-1}\beta\gamma(\xi^{-2}(h_{2})1_{3})n_{(1)})\\
&&\cdot (\gamma S^{-1}(1_{1})\cdot \mu^{-1}(m))\\
&=& \nu(\gamma^{-1}\beta\gamma \xi^{-1}(h_{1})\cdot (\gamma^{-1}\beta\gamma(1_{2})\cdot \mu^{-1}(n_{(0)})))\otimes (\beta^{-1}\delta\gamma^{-1}\beta\gamma\xi^{-1}(h_{2})\\
&&(\beta^{-1}\delta\gamma^{-1}\beta\gamma(1_{3})\beta^{-1}\xi^{-1}(n_{(1)}))\cdot (\gamma S^{-1}(1_{1})\cdot \mu^{-1}(m))\\
&=&\nu(\gamma^{-1}\beta\gamma\xi^{-1}(h_{1})\cdot (1_{2}\cdot \nu^{-2}(n_{(0)})))\otimes (\beta^{-1}\delta\gamma^{-1}\beta\gamma\xi^{-1}(h_{2})\\
&&\beta^{-1}((\delta(1_{3})\xi^{-2}(n_{(1)}))\gamma S^{-1}(1_{1}))\cdot m\\
&=& \gamma^{-1}\beta\gamma(h_{1})\cdot \nu(n_{(0)})\otimes (\beta^{-1}\delta\gamma^{-1}\beta\gamma\xi^{-1}(h_{2})\beta^{-1}(n_{(1)}))\cdot m\\
&=& h\cdot c_{M\otimes N}(m\otimes n).
\end{eqnarray*}

Secondly, we check that $c_{M, N}$ is an $H$-comodule map as
 follows:
\begin{eqnarray*}
&&\rho_{N\otimes M}\ci c_{M,N}(m\otimes n) \\
&=&((\nu(n_{(0)}))_{<0>}\otimes( \beta^{-1}(n_{(1)})\c\mu^{-1}(m))_{(0)})
\otimes( \beta^{-1}(n_{(1)})\c\mu^{-1}(m))_{(1)}\\
&&(\nu(n_{(0)}))_{<1>}\\
&=&(n_{(0)}\otimes \beta^{-1}\xi(n_{(1)21})\c\mu^{-1}(m_{(0)}))
\otimes(\xi(n_{(1)22})\xi^{-1}(m_{(1)}))\\
&&(\alpha\beta^{-1}S^{-1}\xi(n_{(1)12})\alpha\beta^{-1}\xi(n_{(1)11}))\\
&=&\varepsilon(n_{(1)11}1'_{1}1_{2})(n_{(0)}\otimes \beta^{-1}(n_{(1)12}1'_{2})\cdot\mu^{-1}(m_{(0)}))\\
&&\otimes(n_{(1)2}\xi^{-1}(m_{(1)}))\alpha\beta^{-1}S^{-1}(1_{1})\\
&=&(n_{(0)}\otimes \beta^{-1}\xi^{-1}(n_{(1)1}1_{2})\cdot\mu^{-1}(m_{(0)}))\\
&&\otimes(n_{(1)2}\xi^{-1}(m_{(1)}))\alpha\beta^{-1}S^{-1}(1_{1})\\
&=& (n_{(0)}\otimes \beta^{-1}(n_{(1)1})\cdot \mu^{-1}(m_{(0)})\otimes \xi(n_{(1)2})m_{(1)}\\
&=&(\nu( n_{(0)(0)})\otimes \beta^{-1}(n_{(0)(1)})\c\mu^{-1}(m_{(0)}))
\otimes n_{(1)}m_{(1)}\\
&=&(c_{M,N}\otimes id)((m_{(0)}\otimes n_{(0)})\otimes n_{(1)}m_{(1)})
=(c_{M,N}\otimes id)\rho(m\otimes n).
\end{eqnarray*}

Finally we will check Eqs.(4.4) and (4.5). On the one hand,
\begin{eqnarray*}
&&a^{-1}_{{}^{M\otimes N}P, M, N}\circ c _{M\otimes N, P}\circ a^{-1}_{M, N, P}(m\otimes (n\otimes p))\\
&=&a^{-1}_{{}^{M\otimes N}P, M, N}\circ c _{M\otimes N, P}((\mu^{-1}(m)\otimes n)\otimes \varsigma(p))\\
&=&a^{-1}_{{}^{M\otimes N}P, M, N}(\varsigma^{2} ( p_{(0)} )\otimes \gamma^{-1}\beta^{-1}\gamma\delta^{-1}\xi( p_{(1)} )
\c(\mu^{-2}(m)\otimes \nu^{-1}(n)))\\
&=&(c _{M, {}^NP}\otimes id_N)((\mu^{-1}(m)\otimes \varsigma(p_{(0)}))
\otimes \delta^{-1}\xi(p_{(1)})\c n)\\
&=&(c _{M, {}^NP}\otimes id_N)\circ a^{-1}_{M, {}^NP, N}\circ (id _M\otimes c _{N, P})(m\otimes(n\otimes p)).
\end{eqnarray*}
On the other hand,
 \begin{eqnarray*}
 &&a_{{}^MN, {}^MP, M} \circ  c _{M, N\otimes P}\circ a_{M, N, P}((m\otimes n)\otimes p)\\
 &=&a_{{}^MN, {}^MP, M} \circ  c _{M, N\otimes P}(\mu(m)\otimes( n\otimes\varsigma^{-1} (p)))\\
 &=&\nu^{2}(n_{(0)})\otimes \varsigma(\varsigma^{-1}(p)_{(0)})
 \otimes (\beta^{-1}(\varsigma^{-1}(p)_{(1)})\c(\beta^{-1}\xi^{-1} (n_{(1)})\c \xi^{-1}(m)))\\
 &=&(id _{{}^MN }\otimes c _{M, P})
 ((\nu^{2}(n_{(0)})\otimes \beta^{-1}(n_{(1)})\c \mu^{-1}(m))\otimes p)\\
 &=&(id _{{}^MN }\otimes c _{M, P})\circ a_{{}^MN, M, P}\circ (c _{M, N}\otimes id_P)
 ((m\otimes n)\otimes p)).
\end{eqnarray*}

 The proof is completed.
 \hfill $\blacksquare$
\\

{\bf Lemma 4.5.} The map $c_{M,N}$
defined by $c_{M,N}(m\otimes n)=\nu(n_{(0)})
\otimes \beta^{-1}(n_{(1)})\c \mu^{-1}(m)$ is bijective; with inverse
 $${c}_{M,N}^{-1}(n\otimes m)=\beta^{-1}(S(n_{(1)}))\c\mu^{-1}(m) \otimes \nu (n_{(0)}).$$

{\bf Proof.} First, we prove $c_{M,N}{c}_{M,N}^{-1}=id$. For all $m\in M, n\in
  N$, we have
 \begin{eqnarray*}
&&c_{M,N}{c}_{M,N}^{-1}(n \otimes m)\\
&=&c_{M,N}(\beta^{-1}S(n_{(1)})\c\mu^{-1}(m)\otimes\nu(n_{(0)}))\\
&=&\nu(\nu(n_{(0)})_{(0)})\otimes \beta^{-1}(\nu(n_{(0)})_{(1)})\c \mu^{-1}(\beta^{-1}S(n_{(1)})\c\mu^{-1}(m))\\
&=&\nu^{2}(n_{(0)(0)})\otimes \beta^{-1}((n_{(0)(1)})S\xi^{-1}(n_{(1)}))\c\mu^{-1}(m)\\
&=& \alpha\gamma^{-1}\beta^{-1}\gamma\delta^{-1}\beta\alpha^{-1}(1_1)\rhd \nu(n_{<0>})\otimes \alpha^{-1}(1_2\varepsilon_{t}(n_{<1>}))\cdot m\\
&=& 1_1\rhd \nu(n)\otimes \beta^{-1}\delta^{-1}\gamma^{-1}\beta\gamma \alpha^{-1}(1_2)\cdot m\\
&=& n\otimes m.
\end{eqnarray*}

The fact that $ c_{M,N}^{-1} c_{M,N}=id$ is similar.
 This completes the proof. \hfill $\blacksquare$
\\

 Let $H$ be a weak monoidal Hom-Hopf algebra and
 $G={\sl Aut}_{wmHH}(H)
 \times {\sl Aut}_{wmHH}(H)$.
 Define $\mathcal {WMHYD}(H)$ as the
  disjoint union of all $_{H}\mathcal {WMHYD}^{H}(\alpha,\beta)$
 with $(\alpha,\beta)\in G$. If we endow $\mathcal {WMHYD}(H)$
 with tensor product shown in Proposition 4.1,
 then $\mathcal {WMHYD}(H)$ becomes
 a monoidal category with unit $H_{t}$.

 Define a group homomorphism
 $\,\,\varphi: G\rightarrow Aut(\mathcal {WMHYD}(H)),
 \,\,\,\,\,\,\,\,
 (\alpha, \beta) \,\,\mapsto \,\,\varphi(\alpha,\beta)\,\,$
 on components as follows:
\begin{eqnarray*}
\varphi_{(\alpha,\beta)}: {_{H}}\mathcal {WMHYD}^{H}(\gamma,\delta)&\rightarrow&
{_{H}}\mathcal {WMHYD}^{H}((\alpha,\beta)\ast
(\gamma,\delta)\ast (\alpha,\beta)^{-1}),\\
\quad \quad \quad \quad \quad
 \quad \varphi_{(\alpha,\beta)}(N)&=& {}^{(\alpha,\beta)} N,
\end{eqnarray*}
and the functor $\varphi_{(\alpha,\beta)}$ acts as identity on morphisms.\\

 The braiding in $\mathcal {WMHYD}(H)$
 is given by the family $\{c_{M,N}\}$
 in Proposition 4.4.
 So we get the following main theorem of this article.
\\

{\bf Theorem 4.6.} $\mathcal {WMHYD}(H)$
is a braided $T$-category over $G$.
\\

\section*{ACKNOWLEDGEMENTS}
 The work was partially supported by  the Fundamental Research Funds for the Central Universities  (NO. 3207013906),
 and the NSF of China (NO. 11371088), and the NSF of Jiangsu Province (NO. BK2012736).
  \\

\end{document}